\documentclass[onefignum,onetabnum]{siamart171218}

\usepackage{lipsum}
\usepackage{booktabs}
\usepackage{amsfonts}
\usepackage{graphicx}
\usepackage{epstopdf}
\usepackage{subcaption}
\usepackage{algorithmic}
\usepackage{amsmath}
\usepackage[a4paper]{geometry}

\DeclareMathOperator*{\argmin}{arg\,min}

\ifpdf
  \DeclareGraphicsExtensions{.eps,.pdf,.png,.jpg}
\else
  \DeclareGraphicsExtensions{.eps}
\fi

\usepackage{todonotes}
\usepackage{changes}

\newsiamremark{remark}{Remark}
\newsiamremark{hypothesis}{Hypothesis}
\crefname{hypothesis}{Hypothesis}{Hypotheses}
\newsiamthm{claim}{Claim}

\headers{Practical shift choice in the shift-and-invert method}{A. Katrutsa, M. Botchev, and I. Oseledets}

\title{Practical shift choice in the shift-and-invert Krylov subspace evaluations of the matrix exponential\thanks{Submitted to the editors DATE.
\funding{The first and the third authors are supported by RFBR grant 18-31-20069 mol\_a\_ved. The work of the second author is supported by Russian Science Foundation grant No.~19-11-00338. }}}

\author{Alexandr Katrutsa\footnotemark[2]
\and Mike Botchev\footnotemark[3]~\footnotemark[4]
\and Ivan Oseledets\footnotemark[2]~\footnotemark[4]
\footnotetext[2]{S\MakeLowercase{kolkovo} I\MakeLowercase{nstitute of} S\MakeLowercase{cience and }T\MakeLowercase{echnology}, M\MakeLowercase{oscow}, R\MakeLowercase{ussia (\email{aleksandr.katrutsa@phystech.edu}, \email{i.oseledets@skoltech.ru}).}}
\footnotetext[3]{K\MakeLowercase{eldysh} I\MakeLowercase{nstitute of} A\MakeLowercase{pplied} M\MakeLowercase{athematics}, R\MakeLowercase{ussian} A\MakeLowercase{cademy of }S\MakeLowercase{ciences}, M\MakeLowercase{oscow}, R\MakeLowercase{ussia 
  (\email{botchev@ya.ru}).}}
  \footnotetext[4]{M\MakeLowercase{archuk} I\MakeLowercase{nstitute of }N\MakeLowercase{umerical} M\MakeLowercase{athematics}, R\MakeLowercase{ussian} A\MakeLowercase{cademy of} S\MakeLowercase{ciences}, M\MakeLowercase{oscow}, R\MakeLowercase{ussia}}
}

\usepackage{amsopn}

\makeatletter
\newcommand*{\addFileDependency}[1]{% argument=file name and extension
  \typeout{(#1)}% latexmk will find this if $recorder=0 (however, in that case, it will ignore #1 if it is a .aux or .pdf file etc and it exists! if it doesn't exist, it will appear in the list of dependents regardless)
  \@addtofilelist{#1}% if you want it to appear in \listfiles, not really necessary and latexmk doesn't use this
  \IfFileExists{#1}{}{\typeout{No file #1.}}% latexmk will find this message if #1 doesn't exist (yet)
}
\makeatother

\newcommand{\Kk}{\mathcal{K}}
\newcommand{\Rr}{\mathbb{R}}
\newcommand{\Rrnn}{\Rr^{n\times n}}

\begin{document}

\maketitle

% REQUIRED
\begin{abstract}
    We propose two methods to find a proper shift parameter in the shift-and-invert method for computing matrix exponential matrix-vector products.
    These methods are useful in the case of matrix exponential action has to be computed for a number of vectors.
    The first method is based on the zero-order optimization of the mean residual norm for a given number of initial vectors. 
    The second method processes the vectors one-by-one and estimates, for each vector, the derivative of the residual norm as a function of the shift parameter.
    The estimated derivative value is then used to update the shift value for the next vector.
    To demonstrate the performance of the proposed methods we perform numerical experiments for two-dimensional non-stationary convection-diffusion equation with discontinuous coefficients and two-dimensional anisotropic diffusion equation.
    % The total running times for the range of dimensions are provided.
    % The experiments show that extra costs for optimization or derivative estimation pay off already for a small number of vectors to be processed.
\end{abstract}

% REQUIRED
\begin{keywords}
  matrix exponential, shift-and-invert method, rational Krylov method, zero-order optimization, Brent method 
\end{keywords}

% REQUIRED
\begin{AMS}
  65F30, 65F60, 65F10, 65N22, 65L05
\end{AMS}

\section{Introduction}
\label{sec::intro}
Computation of the matrix exponential actions on vectors for large and sparse matrices is an important task occurring, e.g., in time integration of large dynamical systems~\cite{al2011computing,gajic2003linear,gallivan1994pade}, network analysis~\cite{de2019analysis,benzi2015limiting,gilson2018framework}, Markov chain modeling~\cite{metzler2018linear,bladt2017matrix,sidje1999numerical} and many other problems.
Rational Krylov subspace methods~\cite{guttel2013rational,berljafa2015generalized,druskin2011adaptive,knizhnerman2009optimal} and, in particular, \emph{the shift-and-invert Krylov (SAI Krylov) subspace method}~\cite{moret2004rd,van2006preconditioning}, 
form an efficient class of methods often used for this purpose~\cite{lopez2006analysis,druskin2009solution,druskin2010adaptive}. 
Fast convergence and robust behavior of rational Krylov methods are paid by the necessity to solve linear systems at each Krylov step, sometimes with different matrices.
An attractive property of the SAI Krylov method is that a single pole is involved and, hence, linear systems with just a single matrix have to be solved.
More specifically, assume $A\in\Rrnn$ is a large sparse matrix whose Hermitian part $\frac12(A+A^T)$ is positive semidefinite and we are interested in computing 
\begin{equation}
y(t) = \exp(-tA)v
\label{eq::matexp}
\end{equation}
for given $v\in\Rr^n$ and $t>0$. 
Regular Krylov subspace methods usually employ the Galerkin projection of $y=y(t)$ on the Krylov subspace
\begin{equation}
    \Kk_k( A,v ) =
\mathrm{span} (\, v, A v,  A^2v,
\dots, A^{k-1}v \, ).
\label{eq::krylov_subspace}
\end{equation}
The SAI Krylov method works instead with $\Kk_k( (I+\gamma A)^{-1},v )$, where, $(I+\gamma A)^{-1}$ is the shifted-and-inverted matrix and $\gamma>0$ is a parameter, called a shift or pole, whose value has to be chosen properly.

This study presents two methods to choose the shift parameter $\gamma$ in the SAI Krylov method.
The presented methods are of practical interest if the matrix exponential action
has to be computed for a number of vectors.
% MB:
We assume that the initial vectors are in some sense similar, e.g.,
they belong to a certain subspace or have
% multiple vectors from 
the same distribution (for instance, generated from the normal distribution with the same mean and covariance matrix). 
Examples of such setting are time integration of reaction--diffusion systems~\cite{HundsdorferVerwer:book}, inverse problem modeling based on numerous direct problem  solutions~\cite{Borcea_Druskin_ea2019}, 
nuclear reactor optimization~\cite{StrakhovskayaFedorenko2000,KurchenkovaLebedev2007}. 
% computing the action of matrix exponential on many initial states in non-stationary diffusion problem 
% and computing Hutchinson trace estimator~\cite{todo} of the matrix exponential.
% The last example is important in the problem of graph centrality computation~\cite{todo} and in the problem of comparing data manifolds through Gromov-Wasserstein distance~\cite{todo}.

Although both of the proposed methods to choose $\gamma$ give a performance gain in the case of multiple initial  vectors~$v$, each of them utilizes these vectors in its own way.
In particular, the first method, which we call \emph{``optimize-and-run''}, 
employs a couple of the given initial vectors to find a proper
optimal shift parameter in a pre-processing manner.
Hence, it requires that all these vectors are available in advance.
After this pre-processing stage one can use the determined shift value to compute the matrix exponential actions on the other vectors from the initial vector set.
% from the same distribution as the vectors used to find an optimal shift.
% of the initial state set.
% The larger number of vectors that we proceed, the larger gain we get compared with constant shift used in~\cite{van2006preconditioning}.
Obviously, the larger the initial vector set, the larger the gain obtained by this ``optimize-and-run'' method.
% Eto seichas nevazhno, kakoe znachenie parametra my ispol'zuem, esli optimizatsii net

The second method, which we call \emph{incremental}, processes the 
initial vectors $v$ one-by-one and, hence, does not require the vectors be available beforehand.
% It updates the shift parameter for the next vector using a derivative estimate of the residual norm (as a function of the shift) computed for the current vector. 
For each $v$, it computes derivative estimate of the residual norm as a function of $\gamma$ and updates $\gamma$ accordingly.
To compute the derivative estimate we use a finite difference approximation and propose a specific modification of the SAI Krylov method.
The computational costs required for the estimate in this modified version of the SAI Krylov method are moderate.
% This modification has an admissible overhead compared with SAI Krylov method without a derivative estimate.

Both proposed methods require additional computations to find an  optimum shift value.
Therefore, in the presented experiments we show how large these extra costs are and for how many initial vectors they are paid off.
As test problems, we consider time integration of a convection-diffusion equation with piecewise constant coefficients and an anisotropic diffusion equation.
We test performance of the proposed methods for different mesh sizes and values of~$t$.
% MB:
% Experiments show that both proposed methods outperform a priori chosen shift parameter~\cite{van2006preconditioning} in the case of computing~\eqref{eq::matexp} for multiple vectors $v$.
% Moreover, extra costs for the optimization step or computing derivative estimate are moderate compared to gain that we get.
% MB:
% Naverno, luchshe skazat', tak:
A question then arises whether the additional costs for the optimization are paid off and for how many initial vectors.
We discuss this question in details.

% Namely, we consider non-stationary diffusion equation for piecewise constant coefficient with large step size, different graphs and data manifolds 

% MB - zdes' ostanovilsya 03.07.2019

\subsection{Related work}

The problem of the shift choice in the SAI Krylov subspace method is considered in~\cite{guttel2013rational,moret2014restarted,knizhnerman2010new}.
However, the methods presented in these papers find multiple shifts for general rational Krylov method 
and do not address the shift choice in the SAI Krylov method.
In~\cite{moret2004rd} the value $0.1t$ is used for the single shift.
The %MB adaptive 
choice of the single shift in the SAI Krylov method is discussed in~\cite{van2006preconditioning}.
This study considers the approximation error of the exponential by the restricted rational approximation, corresponding to the SAI Krylov method.
This approximation error depends on shift and, hence, provides an indication 
a shift value most suitable %MB the shift quality 
to achieve a given tolerance.
%MB This allows to find the best shift value for the required tolerance.
The best shift values for a range of tolerances are also presented~\cite{van2006preconditioning}.
We compare the shifts obtained by the presented methods with numerical approximations from~\cite{van2006preconditioning} and demonstrate the gain in the total CPU time.

\section{Shift-and-invert Krylov subspace method (SAI Krylov)}
\label{sec::shift-invert}

This section briefly introduces the SAI Krylov method~\cite{moret2004rd,van2006preconditioning}, provides its pseudocode and implementation details.
The regular Krylov method to compute matrix exponential produces an orthonormal basis $V_m = [v_1, \ldots, v_m] \in \mathbb{R}^{n \times m}$ of the $m$-th order Krylov subspace ~\eqref{eq::krylov_subspace} and approximates matrix exponential in this subspace.
The orthonormal basis is constructed by the Arnoldi process~\cite{saad2003iterative, van2003iterative} that uses the Arnoldi relation
\begin{equation}
    AV_m = V_{m+1}H_{m+1,m} = V_m H_{m, m} + v_{m+1}h_{m+1,m} e_1^{\top},
    \label{eq::arnoldi}
\end{equation}
where $H_{m+1,m} \in \mathbb{R}^{(m+1)\times m}$ is upper-Hessenberg matrix, $H_{m, m}$ is a matrix with the first $m$ rows of the matrix $H_{m+1,m}$, $e_1 = [1, 0, \ldots, 0]\in \mathbb{R}^n$ is a column vector.  
Also, we use the fact that the matrix exponential $\exp(-tA)$ is a matrix polynomial in $A$~\cite{gantmacher2005applications}.
Therefore, the following approximation holds
\begin{equation}
    y(t) = \exp(-tA)v = \exp(-tA)\beta V_m e_1 \approx V_m\exp(-tH_{m, m})\beta e_1 = y_m.
    \label{eq::krylov_approx}
\end{equation}
where $v_1 = v / \beta$ and $\beta = \|v\|_2$.
The well-known fact about computing matrix exponential is that the small eigenvalues of the matrix $A$ are typically more important than the large ones.
However, regular Krylov subspace methods tend to detect the less important
large eigenvalues first~\cite{saad2003iterative,van2003iterative},
whereas its convergence with respect to the important small eigenvalues
can be very slow.
To address this issue, the SAI Krylov method was proposed.
This method replaces the matrix $A$ with the shifted-and-inverted matrix $(I + \gamma A)^{-1}$, where $\gamma >0$ is a shift, and, hence, the orthonormal basis of the Krylov subspace is computed for this shifted-and-inverted matrix.
Therefore, the SAI Krylov method better approximates the largest eigenvalues of the $(I + \gamma A)^{-1}$ which correspond to the smallest eigenvalues of the matrix $A$.

According to the Arnoldi relation for the matrix $(I + \gamma A)^{-1}$
\[
(I + \gamma A)^{-1}\hat{V}_m = \hat{V}_{m+1} \hat{H}_{m+1,m} = \hat{V}_{m+1} \hat{H}_{m,m} + \hat{v}_{m+1}\hat{h}_{m+1,m} e_1^{\top},
\]
the SAI Krylov method builds the orthonormal basis $\hat{V}_m$ and the upper-Hessenberg matrix $\hat{H}_{m+1,m} \in \mathbb{R}^{(m+1) \times m}$.
Other notations are the same as in~\cref{eq::arnoldi}.
% Note that the matrices $V_{m+1}$ and $\hat{H}_{m+1,m}$ constructed here by the SAI method are different than those in relation~\eqref{eq::arnoldi}, so that we slightly abuse the notation.
% To derive the matrix $H_{m,m}$ used in~\cref{eq::krylov_approx} one has to apply to the matrix $\hat{H}_{m,m}$ inverse of the shift-and-invert transformation:
Since $\hat{H}_{m,m}$ is a projection of $(I + \gamma A)^{-1}$, the matrix
\begin{equation}
    H_{m, m} = \frac{1}{\gamma}(\hat{H}^{-1}_{m, m} - I).
    \label{eq::inverse_sai}
\end{equation}
should be an approximation to the matrix $A$ and is used in the SAI Krylov method instead of the matrix $H_{m,m}$ in~\eqref{eq::krylov_approx}.
Different stopping criteria for this process exist, for instance smallness of the residual norm~\cite{botchev2013residual}.
Now all the basic ingredients of the SAI Krylov method are discussed and we summarize them in~\cref{alg::shift-invert}.

Typically, if the Krylov subspace dimension is not too large, the most costly operations in the SAI Krylov method is multiplication of the matrix $(I + \gamma A)^{-1}$ by vectors.
This is equivalent to solving the linear systems with a fixed matrix $I + \gamma A$ and different right-hand sides.
To perform the matrix--vector multiplications with $(I + \gamma A)^{-1}$ efficiently, one can carry out a sparse $LU$ factorization of the matrix $I + \gamma A$ and then use it to solve many linear systems.
In our implementation of~\cref{alg::shift-invert} the sparse $LU$ factorization is performed by the SuperLU library~\cite{li2005overview} wrapped in SciPy~\cite{jones2014scipy} with the default parameters. 
Also the matrix exponentials of the projected matrix from~\eqref{eq::inverse_sai} are computed with the standard method
from the SciPy library.
This method is proposed in~\cite{al2009new} and uses Pad\'e approximation.

\begin{algorithm}[!htb]
\centering
\caption{Shift-and-invert Krylov method to compute~\cref{eq::matexp}}
\label{alg::shift-invert}
\begin{algorithmic}[1]
\REQUIRE{Matrix $A$, shift parameter $\gamma$, vector $v \in \mathbb{R}^n$, time $t$, tolerance $\epsilon$, maximum  number of iterations $k$}
\ENSURE{Result vector $y_m$, residual norm $\|r\|_{\infty}$}
%   \STATE{$A_{\gamma} = $}, where $I$ is identity matrix
  \STATE{Compute sparse $LU$ factorization of $I + \gamma A$}
  \STATE{Initialize zero matrices $V \in \mathbb{R}^{n \times (k+1)}$ and $\hat{H} \in \mathbb{R}^{(k+1) \times k}$}
  \STATE{$\beta = \|v\|_2$}
  \STATE{$V_{:,1} = v / \beta$, where $V_{:,1}$ is the first column of the matrix $V$}
  \FOR{$i=1,\ldots, k$}
    \STATE{$w = (I + \gamma A)^{-1}V_{:, i}$, where the sparse $LU$ factorization of $I + \gamma A$ is used}
    
    \FOR{$j=1,\ldots,i$}
        \STATE{$\hat{h}_{ij} = w^{\top}V_{:, j}$}
        \STATE{$w = w - \hat{h}_{ij}V_{:,j}$}
    \ENDFOR
    
    \STATE{$\hat{h}_{(i+1),i} = \|w\|_2$}
    
    \STATE{Initialize zero vector $e \in \mathbb{R}^{i}$ and assign $e_1 = 1$}
    
    \STATE{Initialize zero vector $d \in \mathbb{R}^i$ and assign $d_{i} = 1$}
    
    \STATE{Compute inverse $\tilde{H} = \hat{H}^{-1}_{1:i, 1:i}$}
    \STATE{$H = \frac{1}{\gamma}\left(\tilde{H} - I\right)$, where $I$ is the $i \times i$ identity matrix}
    
    \STATE{$c = \|(I + \gamma A)w\|_2$}
    \STATE{Create a vector $s = [t/3, 2t/3, t]$}
    \STATE{Initialize zero vector $r \in \mathbb{R}^3$}
    \FOR{$j=1,\ldots,3$} 
        \STATE{$u = \exp(-s_j H)e$}
        \STATE{$r_j = \frac{c}{\gamma} d^{\top}\tilde{H}u$}
    \ENDFOR
    
    \IF{$\|r\|_{\infty} < \epsilon$}
        \STATE{\textbf{break}}
    \ENDIF
    \IF{$i=k$}
        \STATE{\textbf{print}(Number of iterations is exceeded!)}
    \ENDIF
    
    \STATE{$V_{:,(i+1)} = w / \hat{h}_{(i+1),i}$}
  
  \ENDFOR
  \STATE{The dimension of the constructed Krylov subspace $m = i$}
  \STATE{$y_m = \beta V_{:,1:m}v$}
  \RETURN {$y_m$, $\|r\|_{\infty}$}
\end{algorithmic}
\end{algorithm}

\section{Optimization of the shift parameter}
\label{sec::optimization}

% The key ingredient of the SAI Krylov method presented in the previous section is the choice of the shift parameter $\gamma$. 
The proper choice of the shift $\gamma$ reduces the number of Arnoldi iterations in the SAI Krylov method.
In this section, we provide two methods to choose a proper value of $\gamma$.
Both of these methods are of practical interest if one has to solve problem~\cref{eq::matexp} for a number of vectors $v^{(i)}$, $i=1,\ldots,M$.  
% MB:
In this case the additional costs for the optimization are paid off
and a gain with respect to the SAI Krylov method with an non-optimized
reasonable $\gamma$ is achieved.
The minimum number of vectors $M_{\min}$ required to get the gain in the total CPU time for the considered applications is estimated in~\cref{sec::experiments}.
% MB: ostanovilsya zdes' 04.07.2019
The two methods are complementary to each other in the following sense.
The first method is applicable if the set, which the vectors $v^{(i)}$ belong to, is known beforehand.
% It means that we can generate vectors from this set or we already have some number of vectors.
Examples of such sets are the set of vectors generated from the Gaussian distribution with a fixed covariance matrix and the set of random fields generated with the fixed covariance function.
We further refer this method as ``optimize-and-run''.
The second method, which we call the  incremental method, works in an ``on fly'' setting, i.e., the vectors $v^{(i)}$, $i=1,\ldots,M$ do not have to be available beforehand and become available one by one.
It may not be known how are they generated, though we assume that they come from the same set.
% Therefore, it is unknown how are they generated, but we assume that they come from the same set.
In the following subsections, we describe the presented methods and discuss their pros and cons.

\begin{remark}
\label{rem::delta}
Since $\gamma$ usually depends on the time $t$, it is convenient for us
to optimize 
\begin{equation}
    \delta = \gamma / t
    \label{eq::delta_gamma}
\end{equation}
instead of $\gamma$. 
We assume that an optimum value $\delta^*$ (corresponding to an optimum $\gamma^*$) lies in some search interval~$[a, b]$.  
Since usually it is suggested to choose $\gamma$ in the range 
$[0.05t,0.2t]$ (see~\cite{van2006preconditioning}), in all the tests
we take $a=0.01$ and set $b$ to a reasonable known value taken 
from~\cite{van2006preconditioning}.
In the experiments section, the value of $b$ is reported  
for each test problem.
\end{remark}

\subsection{``Optimize-and-run'' method}
\label{sec::opt_run}

The idea of this method is related to the idea from the study~\cite{katrutsa2018optimize}. 
In particular, a stochastic approach is presented to optimize preconditioners for the CG method.
The key idea of this approach is the minimization of a certain objective function with respect to the preconditioner parameter.
This objective function is the error norm after $K$ iterations of the preconditioned CG method averaged over some number of initial guess vectors.
In present work, we propose a method to choose an optimal shift parameter~$\gamma^* = \delta^* t$ in the same manner.
Namely, we propose to solve a similar optimization problem for the parameter $\delta$, but to use the averaged residual norm~\cite{botchev2013residual} as the objective function instead of the averaged error norm.

Formally, we solve the following optimization problem:
\begin{equation}
    \delta^* = \argmin_{\delta \in [a, b]} \frac{1}{N}\sum_{i=1}^N \|r(\delta \; | \; A, x^{(i)}, K, \epsilon, t)\|_{\infty},
    \label{eq::opt_problem}
\end{equation}
where $\|r(\delta \; | \; A, x^{(i)}, K, \epsilon, t)\|_{\infty}$ is the residual norm after performing $K$ Arnoldi iterations of the SAI Krylov method, see~\cref{alg::shift-invert}. 
We consider the residual norm as the function of parameter $\delta$ and minimize average of the residual norms for $N$ trial vectors $x^{(1)}, \ldots x^{(N)}$ with respect to $\delta$.
According to our assumption, the feasible set of $\delta$ is the interval $[a, b]$.
The other parameters like matrix~$A$, time~$t$, tolerance $\epsilon$ and number of Arnoldi iterations $K$ are fixed.

Since the target variable $\delta$ is scalar, a sufficiently good approximation of the minimum can be obtained by a zero-order optimization method~\cite{conn2009introduction}.
Every iteration of such methods requires only the value of the objective function in the current point.
To compute objective function in the point $\delta$, we use~\cref{alg::shift-invert} with corresponding parameters $A$, $t$, $x^{(i)}$, $\epsilon$ and according to the equation~\cref{eq::delta_gamma} use shift $\gamma = \delta t$.
To solve problem~\cref{eq::opt_problem}, we choose Brent method~\cite{brent2013algorithms}, which combines inverse quadratic approximation of the objective function and bisection strategies.
Its computational cost is discussed further.

Let $C_{LU}$ be the cost of the single sparse $LU$ factorization and $C_A$ be the cost of Arnoldi iteration.
In fact, the costs of one Arnoldi iteration grow with the iteration number and by $C_A$ we mean the average costs per iteration, so that the costs for $K$ Arnoldi iterations are $KC_A$.
In our case, the most costly operation in every iteration of Brent method is computing objective function. 
At the same time, the most costly subroutines in computing objective function (see~\cref{alg::shift-invert}) are single sparse $LU$ factorization of the matrix $A + \gamma I$ (line 3 in~\cref{alg::shift-invert}) and  $NK$ Arnoldi iterations (lines 7--27 in~\cref{alg::shift-invert}).
Then the total cost of the optimization procedure is 
\[
(C_{LU} + NKC_A)s,
\]
where $s$ is the number of iterations performed by Brent method to achieve the required tolerance.

% Let $C_{LU}$ be the cost of the single sparse $LU$ factorization and $C_A$ be the cost of one Arnoldi iteration.
% Here, we assume that the most costly part of Arnoldi iteration is solving linear system with matrix $A + \gamma I$. 
% To pay off costs for optimization, we have to get optimum shift $\gamma^*$ such that the costs for processing $M$ vectors $v^{(1)}, \ldots, v^{(M)}$ plus the cost of the optimization are less than the costs of the processing the same $M$ vectors but using the other a priori chosen ``reasonable'' shift~$\bar{\gamma}$.
To pay off the optimization costs, we have to get an optimum shift $\gamma^*$ that significantly reduces costs for processing $M$ vectors $v^{(1)}, \ldots, v^{(M)}$.
In particular, the costs for optimization stage and further processing $M$ vectors with the optimum shift $\gamma^*$ have to be less than the costs for the processing of the same $M$ vectors but using a priori chosen, non-optimized, shift~$\bar{\gamma}$.
By processing we mean here computing matrix exponential action with the desired tolerance $\epsilon$.
The proper choice of the shift can significantly reduce the number of Arnoldi iterations for processing the vectors $v^{(1)},\ldots,v^{(M)}$.
Since these vectors come from the same set, we assume that the number of Arnoldi iterations for every vector $v^{(i)}$ is almost the same.
Therefore, denote by $\bar{K}$ and $K^*$ the number of Arnoldi iterations corresponding to the shifts $\bar{\gamma}$ and $\gamma^*$ respectively. 
Thus, the minimum number of vectors $M_{\min}$ is the smallest number such that 
\[
(s + 1)C_{LU} + (NK + M_{\min}K^*)C_A < C_{LU} + C_A M_{\min}\bar{K}.
\]
Here the right-hand side is the cost of the SAI Krylov method used some fixed a priori chosen shift $\bar{\gamma}$ and the left-hand side is the total cost of the ``optimize-and-run'' method including optimization stage.
We can observe that $M_{\min}$ is small if the two conditions hold:
\begin{enumerate}
    \item[1)] the matrix $A$ has specific sparsity pattern that makes the costs for sparse $LU$ factorization affordable
    \item[2)] the optimum shift $\gamma^*$ leads to a significant reduction of the number of Arnoldi iterations.
\end{enumerate}
Section~\ref{sec::experiments} presents the test problems where these conditions hold.

At each iteration of the Brent method, one sparse $LU$ factorization has to be computed.
In \cref{sec::experiments} the costs of the SAI Krylov method including the optimization stage (the ``optimize-and-run'' method) are compared with the costs of the non-optimized SAI Krylov method. 
In \cref{sec::experiments} we also discuss why the additional costs required to find $\delta^*$ are moderate compared with the obtained gain.
The other important parameter in optimization procedure is the number of Arnoldi iterations $K$. 
This parameter should be chosen as smallest as possible but beyond a convergence stagnation phase typically observed at first iterations. 
% This parameter is chosen in a way to be as smallest as possible and to provide residual norm decreasing simultaneously.
In particular, if $K$ is too small, optimization procedure can not find $\delta^*$ that provides a faster convergence, but if $K$ is too large, then optimization procedure becomes too costly.
In~\cref{sec::experiments} we report the values of $K$ for the considered test problems and show how the choice of $K$ affects the optimization costs and obtained speed up.
To highlight the main steps of the ``optimize-and-run'' method we list them in~\cref{alg::opt_run}.

\begin{algorithm}[!htb]
\centering
\caption{``Optimize-and-run'' method}
\label{alg::opt_run}
\begin{algorithmic}[1]
\STATE{Generate $N$ trial vectors or take $N$ given ones}
\STATE{Fix the number of iterations $K$ and tolerance $\epsilon$}
\STATE{Solve problem~\cref{eq::opt_problem} with $N$ trial vectors and get optimal $\delta^*$}
\STATE{Use $\gamma^* = \delta^*t$ to compute~\cref{eq::matexp} for any other vectors from the same set as the trial vectors from the line 1}
\end{algorithmic}
\end{algorithm}

% \begin{enumerate}
%     \item Generate $N$ trial vectors or take $N$ given ones
%     \item Fix the number of iterations $K$ and tolerance $\epsilon$
%     \item Solve problem~\cref{eq::opt_problem} with $N$ trial vectors and get optimal $\delta^*$
%     \item Use $\gamma^* = \delta^*t$ to compute~\cref{eq::matexp} for any other vectors from the same set as the trial vectors from the step 1.
% \end{enumerate}

\subsection{Incremental method}
Assume we apply the SAI Krylov method to the vector $v^{(i)}$ and use the shift $\gamma^{(i)}$.
The main idea of the incremental approach is to update the shift parameter $\gamma^{(i)}$, while computing $\exp(-tA)v^{(i)}$, so that the updated shift $\gamma^{(i+1)}$ leads to a faster convergence of the SAI Krylov method for the next vector $v^{(i+1)}$.
We expect that this method incrementally speed up convergence of the SAI Krylov method for every next vector.
To update $\gamma^{(i)}$ we can estimate the derivative~$r^{(i)}_{\gamma}$ of the residual norm as a function of $\gamma^{(i)}$.
This estimation requires additional costs which should be paid off with the obtained convergence speed up.
In~\cref{sec::experiments} we provide total CPU time comparison and show for which number of vectors this method becomes attractive.

The incremental method gives a gain in the total CPU time for processing $M$ vectors with respect to the non-optimized SAI Krylov method if
\begin{equation}
  MC_{LU} + (M-1)C_{r_{\gamma}} + (\tilde{K}_1 + \ldots \tilde{K}_M)C_A < C_{LU} + M\bar{K},
  \label{eq::incremental_cost}
\end{equation}
where the left-hand side is the total costs of the incremental method and the right-hand side is the cost of the SAI Krylov method, $C_{r_{\gamma}}$ is the cost to estimate the derivative $r_{\gamma}$ and $\tilde{K}_i$ is the number of Arnoldi iterations to process the $i$-th vector by the incremental method.
Since the shift parameter is non-constant, the number of Arnoldi iterations also varies with the vectors $v^{(i)}$.
Inequality~\cref{eq::incremental_cost} holds if a significant reduction of the number of Arnoldi iterations is obtained with the incremental updates.  

To estimate $r^{(i)}_{\gamma}$ we should have liked to use automatic differentiation tools, like PyTorch~\cite{paszke2017automatic}, Autograd~\cite{maclaurin2015autograd}, etc, but unfortunately these tools do not support necessary operations with sparse matrices, yet.
Therefore, we use a finite difference approximation approach, where a sparse $LU$ factorization is computed only once to decrease additional costs.
In fact, one preconditioned Richardson iteration is used to solve the SAI system for $\gamma + \Delta \gamma$, where the preconditioner is $I + \gamma A$ (lines 10--12 in~\cref{alg::der_eval}).
According to~\cite{kelley1995iterative}, we use $\Delta \gamma = 10^{-7}$ to compute the finite difference approximation.
A detailed description of the proposed modification of the SAI Krylov method that estimates derivative $r^{(i)}_{\gamma}$ is given as~\cref{alg::der_eval}.

\begin{algorithm}[!htb]
\centering
\caption{Estimate derivative $r_{\gamma}$ of the residual norm}
\label{alg::der_eval}
\begin{algorithmic}[1]
\REQUIRE{Matrix $A$, shift parameter $\gamma$, vector $v \in \mathbb{R}^n$, time $t$, tolerance $\epsilon$, maximum number of iterations $k$}
\ENSURE{Result vector $y_m$, derivative estimation $r_{\gamma}$}
  \STATE{$\gamma' = \gamma + 10^{-7}$}
%   \STATE{$A_{\gamma} = I + \gamma A$}, where $I$ is identity matrix
  \STATE{Compute sparse $LU$ factorization of $I + \gamma A$}
  \STATE{Initialize zero matrices $V, V' \in \mathbb{R}^{n \times (k+1)}$ and $\hat{H}, \hat{H}' \in \mathbb{R}^{(k+1) \times k}$}
%   \STATE{Initialize zero matrices $V' \in \mathbb{R}^{n \times (k+1)}$ and $\hat{H}' \in \mathbb{R}^{(k+1) \times k}$}
  \STATE{$\beta = \|v\|_2$}
  \STATE{$V_{:,1} = v / \beta$, where $V_{:,1}$ is the first column of the matrix $V$}
  \STATE{$V'_{:,1} = v / \beta$, where $V_{:,1}$ is the first column of the matrix $V'$}
  \FOR{$i=1,\ldots, k$}
    \STATE{$w = (I + \gamma A)^{-1}V_{:, i}$, where the  sparse $LU$ factorization of $I + \gamma A$ is used}
    
    \STATE{$w' = (I + \gamma A)^{-1}V'_{:, i}$, where the sparse $LU$ factorization of $I + \gamma A$ is used}
    \STATE{$\bar{r} = V'_{:,i} - (w' + \gamma'Aw')$}
    \STATE{$w' = w' + (I + \gamma A)^{-1}\bar{r}$}
    
    \FOR{$j=1,\ldots,i$}
        \STATE{$\hat{h}_{ij} = w^{\top}V_{:, j}$ and $\hat{h}'_{ij} = w'^{\top}V'_{:, j}$}
        \STATE{$w = w - \hat{h}_{ij}V_{:,j}$ and $w' = w' - \hat{h}'_{ij}V'_{:,j}$}
    \ENDFOR
    
    \STATE{$\hat{h}_{(i+1),i} = \|w\|_2$ and $\hat{h}'_{(i+1),i} = \|w'\|_2$}
    
    \STATE{Initialize zero vector $e \in \mathbb{R}^{i}$ and assign $e_1 = 1$}
    
    \STATE{Initialize zero vector $d \in \mathbb{R}^i$ and assign $d_{i} = 1$}
    
    \STATE{Compute inverses $\tilde{H} = \hat{H}^{-1}_{1:i, 1:i}$ and $\tilde{H}' = \hat{H}'^{-1}_{1:i, 1:i}$}
    \STATE{$H = \frac{1}{\gamma}\left(\tilde{H} - I\right)$ and $H' = \frac{1}{\gamma'}\left(\tilde{H}' - I\right)$, where $I$ is the $i \times i$ identity matrix}
    
    \STATE{$c = \|(I + \gamma A)w\|_2$ and $c' = \|(I + \gamma' A)w\|_2$}
    \STATE{Create a vector $s = [t/3, 2t/3, t ]$}
    \STATE{Initialize zero vector $r \in \mathbb{R}^3$ and $r' \in \mathbb{R}^3$}
    \FOR{$j=1,\ldots,3$} 
        \STATE{$u = \exp(-s_j H)e$ and $u' = \exp(-s_j H')e$}
        \STATE{$r_j = \frac{c}{\gamma} d^{\top}\tilde{H}u$ and $r'_j = \frac{c'}{\gamma'} d^{\top}\tilde{H}'u'$}
    \ENDFOR
    
    \IF{$\|r\|_{\infty} < \epsilon$}
        \STATE{$r_{\gamma} = \frac{\|r'\|_{\infty} - \|r\|_{\infty}}{\gamma' - \gamma}$}
        \STATE{\textbf{break}}
    \ENDIF
    
    \STATE{$V_{:,(i+1)} = w / \hat{h}_{(i+1),i}$ and $V'_{:,(i+1)} = w' / \hat{h}'_{(i+1),i}$}
  
  \ENDFOR
  \STATE{The dimension of the constructed Krylov subspace $m = i$}
  \STATE{$y_m = \beta V_{:,1:m}v$}
  \RETURN {$y_m$, $r_{\gamma}$}
\end{algorithmic}
\end{algorithm}

According to~\cref{rem::delta} we optimize~$\delta = \frac{\gamma}{t}$ instead of $\gamma$.
Therefore, instead of updating~$\gamma^{(i)}$ we update~$\delta^{(i)}$ and use~$\gamma^{(i+1)} = \delta^{(i+1)}t$ for the processing the $(i+1)$-th vector $v^{(i+1)}$. 
To update~$\delta^{(i)}$, we propose~\cref{alg::update_delta} consisting of the following steps: take the midpoint $\delta^{(i)}$ of the given interval~$[\delta^{(i)}_l, \delta^{(i)}_u]$, estimate~$r^{(i)}_{\gamma}$ at point $\gamma^{(i)} = \delta^{(i)}t$, update the search interval taking into account the sign of~$r^{(i)}_{\gamma}$ and compute~$\delta^{(i+1)}$ as the midpoint of the updated interval.
For the first vector~$v^{(1)}$ we use the same search interval as in the ``optimize-and-run'' method~\cref{eq::opt_problem}: $\delta^{(1)}_l = a, \; \delta^{(1)}_u = b$.
Also, if the difference between two sequential $\delta$ is less in modulus than~$10^{-5}$, the incremental method stops updating $\delta$.
This stopping rule reduces the number of sparse $LU$ factorizations and, consequently, the total CPU time.
So, we have two stages in the incremental method. At the first stage the method processes initial vectors and updates $\delta$ according to~\cref{alg::update_delta} until it has converged to some value which we denote $\tilde{\delta}$. 
At the second stage the SAI Krylov method runs with the constant shift $\tilde{\gamma} = \tilde{\delta}t$ for the remaining initial vectors without derivative estimation.

\begin{algorithm}[!htb]
\centering
\caption{Update $\delta$}
\label{alg::update_delta}
\begin{algorithmic}[1]
\REQUIRE{Matrix $A$, left and right bounds of the search interval $\delta^{(i)}_l$ and $\delta^{(i)}_u$, time $t$, tolerance~$\epsilon$, maximum number of iterations $k$, initial vector  $v^{(i)}$}
\ENSURE{Updated left and right bounds of the search interval $\delta^{(i+1)}_l$ and 
$\delta^{(i+1)}_u$}
\STATE{$\delta^{(i)} = \frac{\delta^{(i)}_l + \delta^{(i)}_u}{2}$}
  \STATE {Compute $r^{(i)}_{\gamma}$  with~\cref{alg::der_eval} using $\gamma^{(i)} = \delta^{(i)}t$, $A$, $t$, $\epsilon$, $v^{(i)}$ and $k$}
  \IF{$r^{(i)}_{\gamma} > 0$} 
    \STATE{$\delta^{(i+1)}_u = \delta^{(i)}$}
    \STATE{$\delta^{(i+1)}_l = \delta^{(i)}_l$}
  \ELSE
    \STATE{$\delta^{(i+1)}_l = \delta^{(i)}$}
    \STATE{$\delta^{(i+1)}_u = \delta^{(i)}_u$}
  \ENDIF
  \RETURN {$\delta^{(i+1)}_l, \delta^{(i+1)}_u$}
\end{algorithmic}
\end{algorithm}

The main steps of the incremental method are listed below.
\begin{enumerate}
    \item Take the $i$-th vector $v^{(i)}$. 
    If $i = 1$, then $\delta^{(1)}_l = a, \; \delta^{(1)}_u = b$, else use $\delta^{(i)}_l$ and $\delta^{(i)}_u$ coming from the processing of the $(i-1)$-th vector $v^{(i-1)}$
    \item Compute $\delta^{(i)} = \frac{\delta^{(i)}_l + \delta^{(i)}_u}{2}$ and $\gamma^{(i)} = \delta^{(i)}t$
    \item Solve problem~\cref{eq::matexp} with~\cref{alg::der_eval} and get $r^{(i)}_{\gamma}$
    \item Update $\delta^{(i)}_l$ and  $\delta^{(i)}_u$ with lines 3~-- 9 of~\cref{alg::update_delta}
    \item If $|\delta^{(i)}_l - \delta^{(i)}_u| \leq 10^{-5}$, then stop updating $\delta$ and use $\bar{\delta} = \delta^{(i)}$ for remaining initial vectors.
\end{enumerate}

\section{Numerical experiments}
\label{sec::experiments}

In this section, we test the two proposed methods for choosing the shift in the SAI Krylov method.
The SAI Krylov methods with incorporated strategies to choose a proper shift value adaptively are compared against the SAI Krylov method run with a fixed, a priori chosen reasonable shift value. 
% Namely, the proposed methods of adaptive choice of shift and the choice of some ``reasonable'' fixed shift~\cite{van2006preconditioning} are studied.
% The authors of~\cite{van2006preconditioning} provide the table with approximately optimum shifts for the range of accuracy $\epsilon$ in the SAI Krylov method.
To demonstrate the performance of the proposed methods, we consider a non-stationary convection-diffusion equations with piecewise constant coefficients and a non-stationary anisotropic diffusion equation.
Standard finite difference discretizations of the spatial differential operators in both cases give matrices $A$, whose Hermitian part $\frac{1}{2}(A + A^{\top})$ is positive definite.
The source code can be found at GitHub\footnote{\url{https://github.com/amkatrutsa/shift4sai_krylov}}.

\subsection{Non-stationary convection-diffusion equation with piecewise constant coefficients}
\label{sec::diff_piecewise}
The first problem of interest is time integration of non-stationary convection-diffusion equation with piecewise constant coefficients. 
This initial-boundary value problem reads
\begin{equation}
\begin{aligned}
& \frac{\partial u}{\partial t} = (D_1 u_x)_x + (D_2 u_y)_y + \frac{\mathrm{Pe}}{2} (v_1u_x + v_2u_y + (v_1u)_x + (v_2u)_y),\\
& u(x,y, 0) = u_0, \quad u|_{\partial \Omega} = 0 \quad (x,y)\in\Omega=[0,1]\times[0,1]
\end{aligned}
\label{eq::piecewise_diff}
\end{equation}
where the subscripts $\cdot_{x,y}$ denote the partial derivatives with respect to $x$ and $y$, $u_0 = u_0(x, y)$ is an initial state and $\mathrm{Pe}$ is the P\'eclet number.
We set $\mathrm{Pe} = 1000$ and functions $v_1 (x, y) = x + y$ and $v_2(x, y) = x - y$ in all experiments. 
The coefficients $D_{1,2}$ are discontinuous:
\[
D_1 = \begin{cases}
1000, & \quad (x,y) \in \left[\frac14,\frac34\right]\times \left[\frac14,\frac34\right],
\\
0.1, & \quad \text{otherwise},
\end{cases}
\qquad
D_2 = \frac12 D_1.
\]
% Peclet number $p = 1000$ in all experiments.
The matrix~$A$ is obtained by the standard second order central finite difference approximation, so that the diffusion and convection terms contribute
to the symmetric and skew-symmetric parts of $A$, respectively~\cite{Krukier79}.
Let $n$ be the number of points in both spatial dimensions, then the matrix $A$ is of size~$n^2 \times n^2$. 

To test the proposed methods, we solve the problem~\eqref{eq::piecewise_diff} for different initial states generated from the multivariate normal distribution $\mathcal{N}(\mu, \Sigma)$, where the mean~$\mu$ is a random point from~$\Omega$ and the covariance matrix $\Sigma = 0.05 I$, where $I$ is the $2 \times 2$ identity matrix.
% Therefore, to generate every test initial state~$v^{(i)}, i=1,\ldots,M$ we use some random point from~$\Omega$ as~$\mu$ and fixed covariance matrix~$\Sigma$.
Therefore, to generate every initial vector $v^{(i)}, i=1,\ldots, M$ we use some random point from~$\Omega$ as~$\mu$ and the fixed covariance matrix~$\Sigma$.
In experiments we use~$M=20$ initial vectors.
The tolerance in the SAI Krylov method for this test problem is $\epsilon=10^{-6}$.

% We compare the total running time of the SAI Krylov method using the shift given by proposed methods with the SAI Krylov method using the constant shift $\bar{\gamma} = \bar{\delta}t = 0.1t$  for every given initial state $v^{(i)}$. 
According to~\cite{van2006preconditioning}, the shift $\bar{\gamma} = \bar{\delta}t = 0.1t$ is close to the optimal shift for the chosen tolerance $\epsilon = 10^{-6}$.
Therefore, we compare the performance of the SAI Krylov method with the two adaptive ways to choose the shift and the SAI Krylov method with a priori fixed shift $\bar{\gamma} = 0.1t$.
Also we assume that the optimal $\delta^*$ lies in the interval $[0.01, 0.1]$. 
Therefore, we set $a = 0.01$ and $b = 0.1$ in both proposed methods.
 
\subsubsection{``Optimize-and-run'' method}
\label{sec::opt_run_discont_diff}

In this section we present the comparison of the optimum shift $\gamma^* = \delta^* t$ given by the ``optimize-and-run'' method with the a priori chosen shift $ \bar{\gamma} = \bar{\delta}t = 0.1t$.
We demonstrate that the ``optimize-and-run'' method gives the shift $\gamma^* = \delta^*t$ that provides a smaller total CPU time of the SAI Krylov method than the a priori chosen shift $\bar{\gamma} = 0.1t$.
The total CPU time consists of processing initial vectors and the time spent for optimization.
It turns out that a small number of initial vectors is already enough to gain in the total  CPU time.
We now consider costs of optimization stage in more details and then discuss obtained results.

\paragraph{Optimization cost analysis}
To minimize optimization costs, we use the only one trial initial vector, so $N = 1$.
To be sure that such a choice of $N$ does not affect the optimum $\delta^*$, we study the effect of the different trial initial states used in the optimization stage.
This study shows that even if the optimum $\delta^*$ varies for different trial initial states, the number of iterations $K^*$ obtained for corresponding $\delta^*$ is almost the same.
Therefore, we conclude that using a single initial vector in the optimization stage is meaningful.
Tolerance of Brent method is set to $10^{-5}$ and, therefore, we show only 5 digits of $\delta^*$ in the plots below.

In~\cref{tab::opt_run_par_n200_pe500,tab::opt_run_par_n300_pe1000} we provide the total number of sparse $LU$ factorizations and the total number of Arnoldi iterations that are performed during the optimization stage.
Also, in~\cref{tab::opt_run_par_n200_pe500,tab::opt_run_par_n300_pe1000}  the number of iterations of the SAI Krylov method $K$ to compute the objective function in~\eqref{eq::opt_problem} is given.
If the total number of $LU$ factorizations is $s$, then the total number of Arnoldi iterations is $NKs$, but since we use a single  trial vector, it equals $Ks$.
Note that the larger $t$, the larger number of iterations $K$ is required to compute the objective function in~\cref{eq::opt_problem}.
Below we show how the number of iterations $K$ affects the value of $\delta^*$ and obtained gain. 
For $n=300$ and large~$t$ experiments are quite time consuming, therefore  we present the results only for moderate times $t = 10^{-4}$ and $t = 2 \cdot 10^{-4}$.

\begin{table}[!htb]
    \centering
    \caption{Costs of the ``optimize-and-run'' method to find a proper shift value $\gamma^* = \delta^*t$ for $n = 200$, $A \in \mathbb{R}^{n^2 \times n^2}$,  problem~\eqref{eq::piecewise_diff}}
    \begin{tabular}{ccc}
    \toprule
        & $t = 10^{-4}$ & $t = 4 \cdot 10^{-4}$  \\
        \midrule
    Total number of sparse $LU$ factorizations, $s$ & 17 & 12 \\
    Number of Arnoldi iterations for every initial vector, $K$ & 25 & 70\\
    Total number of Arnoldi iterations & $17 \cdot 25$ & $12 \cdot 70$ \\
    \bottomrule
    \end{tabular}
    \label{tab::opt_run_par_n200_pe500}
\end{table}

\begin{table}[!htb]
    \centering
    \caption{Costs of the ``optimize-and-run'' method to find a proper shift value $\gamma^* = \delta^*t$ for $n = 300$, $A \in \mathbb{R}^{n^2 \times n^2}$, problem~\eqref{eq::piecewise_diff}}
    \begin{tabular}{ccc}
    \toprule
        & $t = 10^{-4}$ & $t = 2 \cdot 10^{-4}$ \\
        \midrule
    Total number of sparse $LU$ factorizations, $s$ & 9 & 18  \\
    Number of Arnoldi iterations for every initial vector, $K$ & 50 & 100 \\
    Total number of Arnoldi iterations  & $9 \cdot 50$ & $18 \cdot 100$ \\
    \bottomrule
    \end{tabular}
    \label{tab::opt_run_par_n300_pe1000}
\end{table}

\paragraph{Performance comparison} 
To show the gain given by the ``optimize-and-run'' method, we report an average number of Arnoldi iterations per initial state and the total CPU time in seconds.
\cref{tab::opt_run_arnoldi_test} shows that the found shift $\gamma^* = \delta^*t$ leads to a significant reduction of the average number of Arnoldi iterations compared with the shift $\bar{\gamma} = 0.1t$.
These average numbers of Arnoldi iterations are almost equal to the $K^*$ and $\bar{K}$, respectively.
The following observation explains such a big difference between $K^*$ and $\bar{K}$.
The piecewise constant coefficients $D_1$ and $D_2$ lead to a stiffness in the matrix $A$.
There is a large discrepancy in the magnitude of the eigenvalues of $A$ caused by the stiffness.
Hence, the ability of Krylov subspace methods to adopt to the discrete character of the spectrum  leads to a strong dependence of the number of Arnoldi iterations on the shift.
Thus, using an optimum shift $\gamma^*$ can significantly reduce the number of Arnoldi iterations.

Now, consider the cost $C_{LU}$ of one sparse $LU$ factorization.
Since matrix $A$ comes from discretization of~\cref{eq::piecewise_diff}, it has a  specific sparsity pattern.
This sparsity leads to a low cost of sparse $LU$ factorization of the matrix $I + \gamma A$.
Therefore, according to our complexity analysis of the ``optimize-and-run'' method presented in~\cref{sec::opt_run}, we can  expect that this method, including optimization, is more efficient than running the SAI Krylov method with a reasonable nonoptimal shift value.

\begin{table}[!ht]
    \centering
    \caption{Average number of Arnoldi iterations over 20 initial vectors to achieve tolerance $10^{-6}$, matrix $A \in \mathbb{R}^{n^2 \times n^2}$, problem~\eqref{eq::piecewise_diff}}
    \begin{tabular}{ccc}
    \toprule
       & $\bar{\delta} = 0.1$ & $\delta^*$  \\
       \midrule
     $n=200, t=10^{-4}$  & $88.4$ & $48.5$ \\
     $n=200, t=4 \cdot 10^{-4}$ & $292.05$  & $93.45$ \\
     $n=300, t=10^{-4}$ &$237.75$  & $100.6$ \\
     $n=300, t=2 \cdot 10^{-4}$ & $517.7$ & $150.75$ \\
     \bottomrule
    \end{tabular}
    \label{tab::opt_run_arnoldi_test}
\end{table}

\begin{figure}[!htb]
    \centering
    \begin{subfigure}[t]{0.48\textwidth}
    \centering
    \includegraphics[scale=0.27]{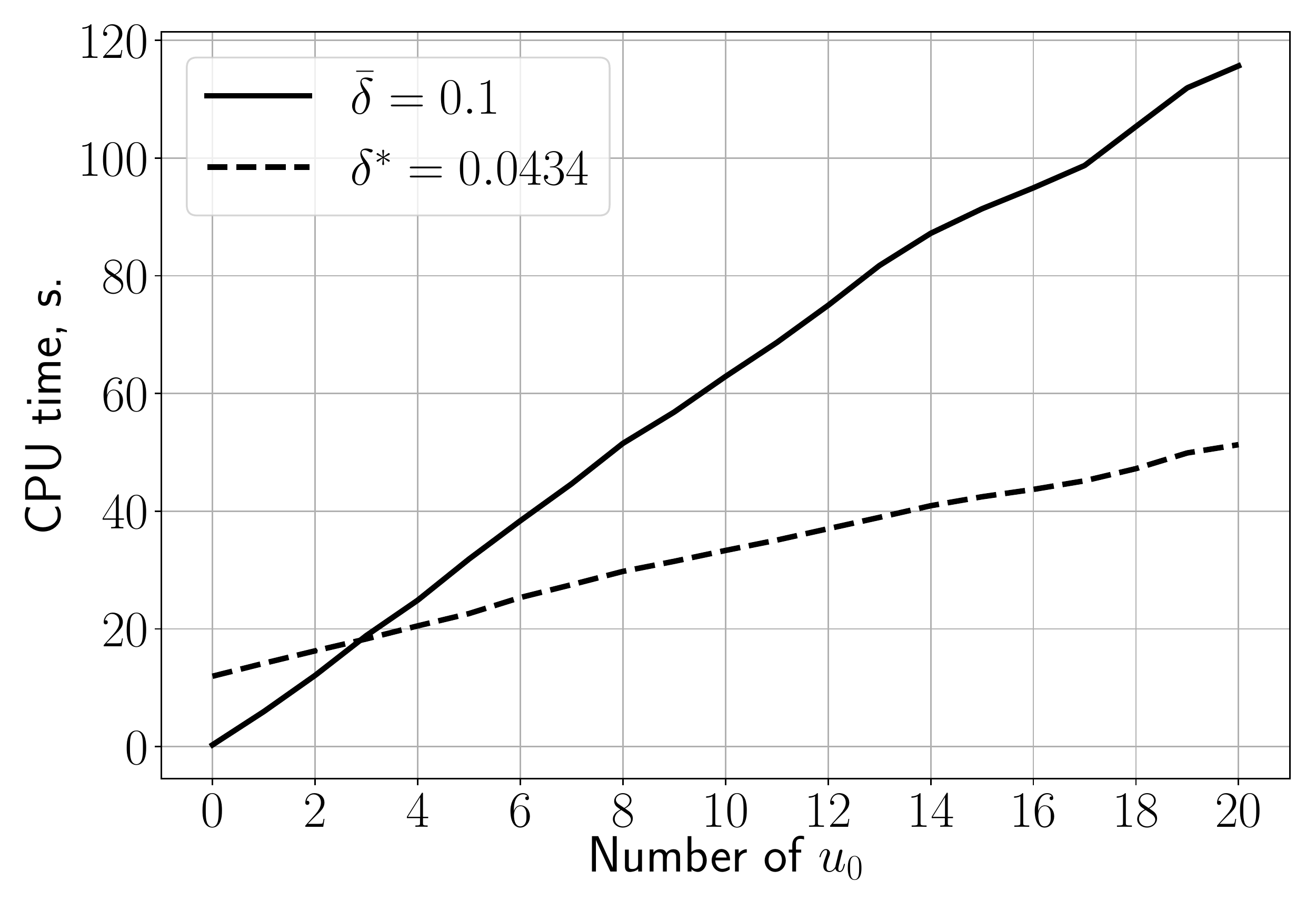}
    \caption{$n = 200, t = 10^{-4}$}
    \label{fig::opt-run_n200_pe500_tau1e-4}
    \end{subfigure}
    ~
    \begin{subfigure}[t]{0.48\textwidth}
    \centering
    \includegraphics[scale=0.27]{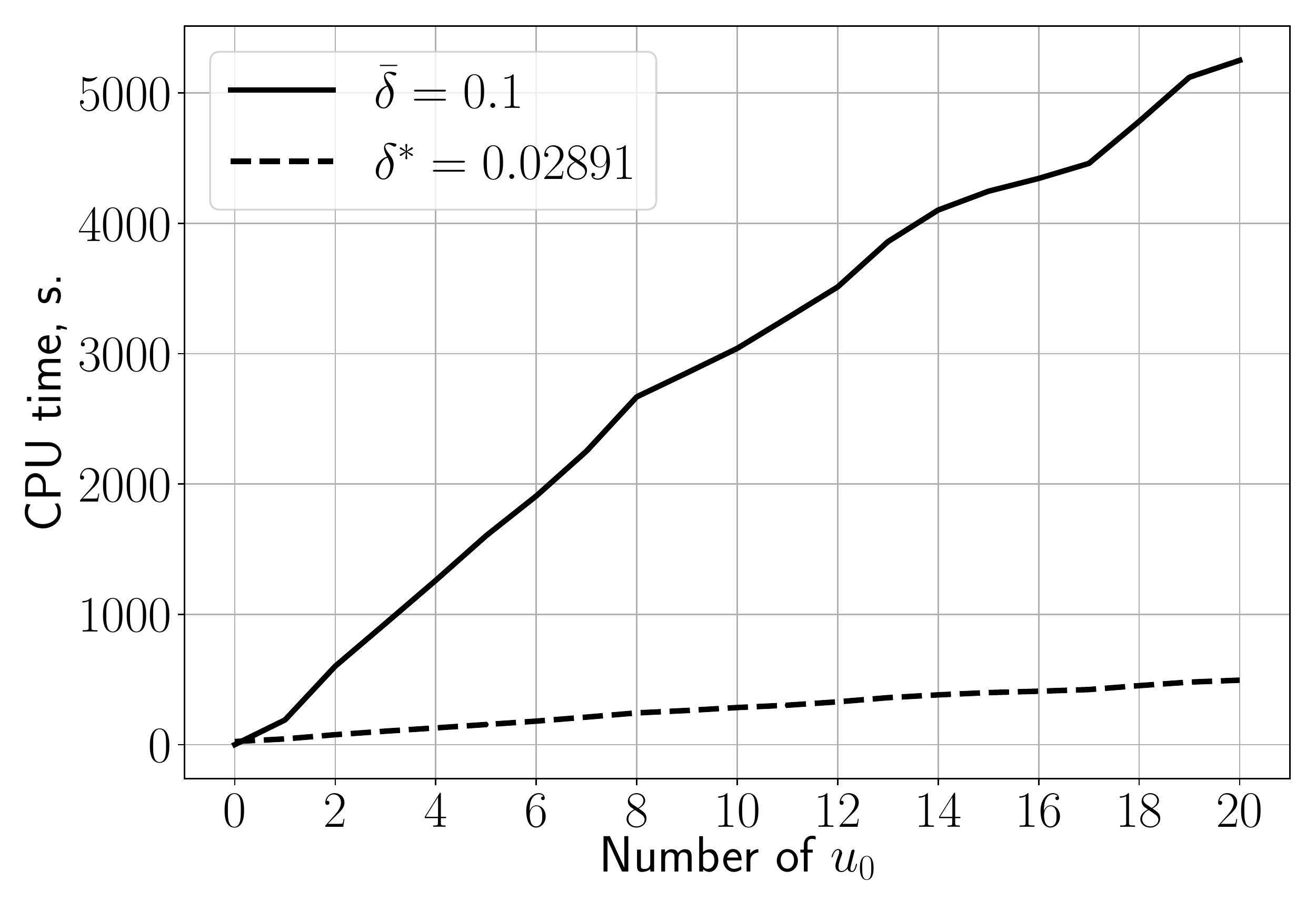}
    \caption{$n = 200, t = 4 \cdot 10^{-4}$}
    \end{subfigure}
    \\
    \begin{subfigure}[t]{0.48\textwidth}
    \centering
    \includegraphics[scale=0.27]{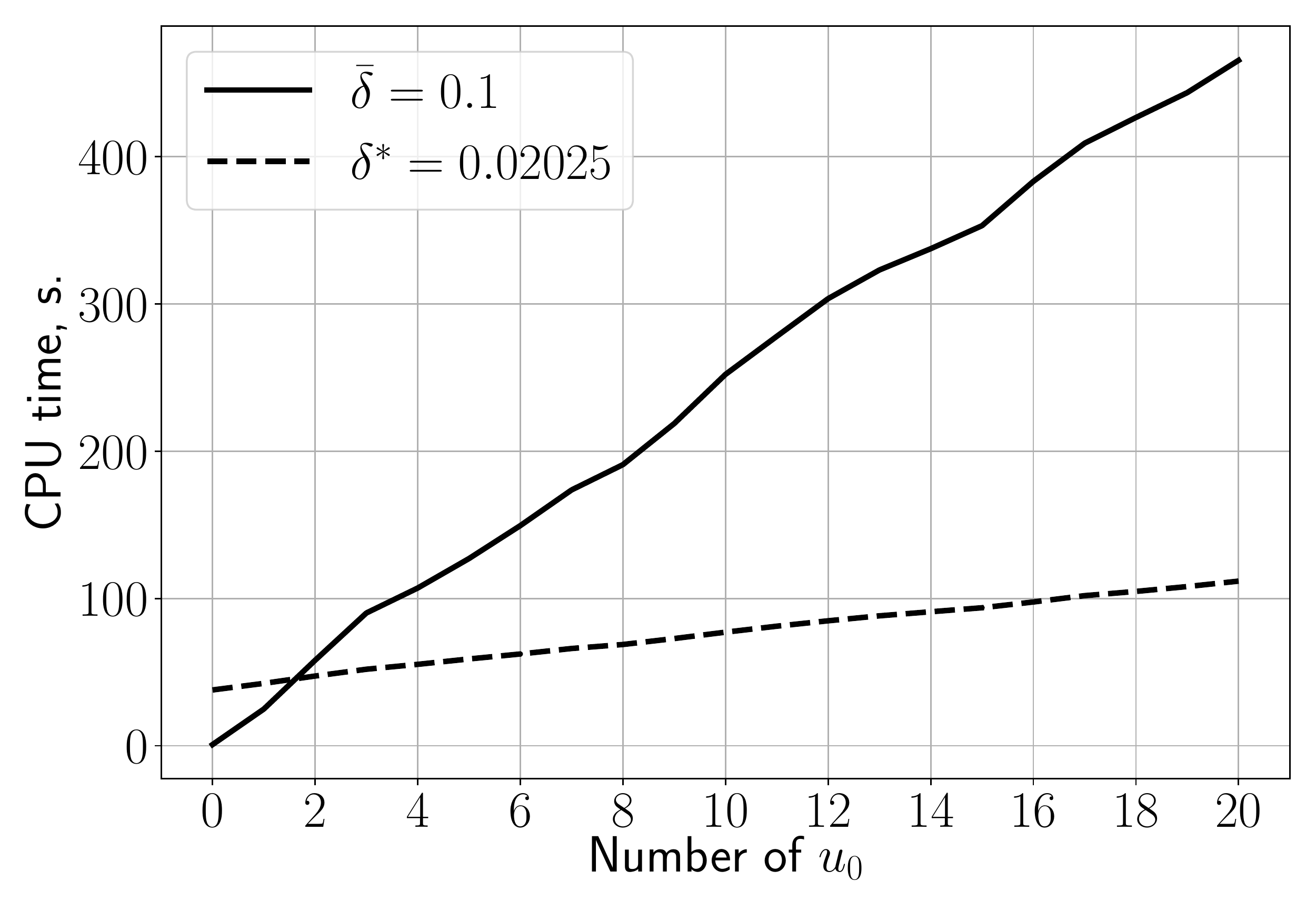}
    \caption{$n = 300, t = 10^{-4}$}
    \label{fig::opt-run_n300_pe1000_tau1e-4}
    \end{subfigure}
    ~
    \begin{subfigure}[t]{0.48\textwidth}
    \centering
    \includegraphics[scale=0.27]{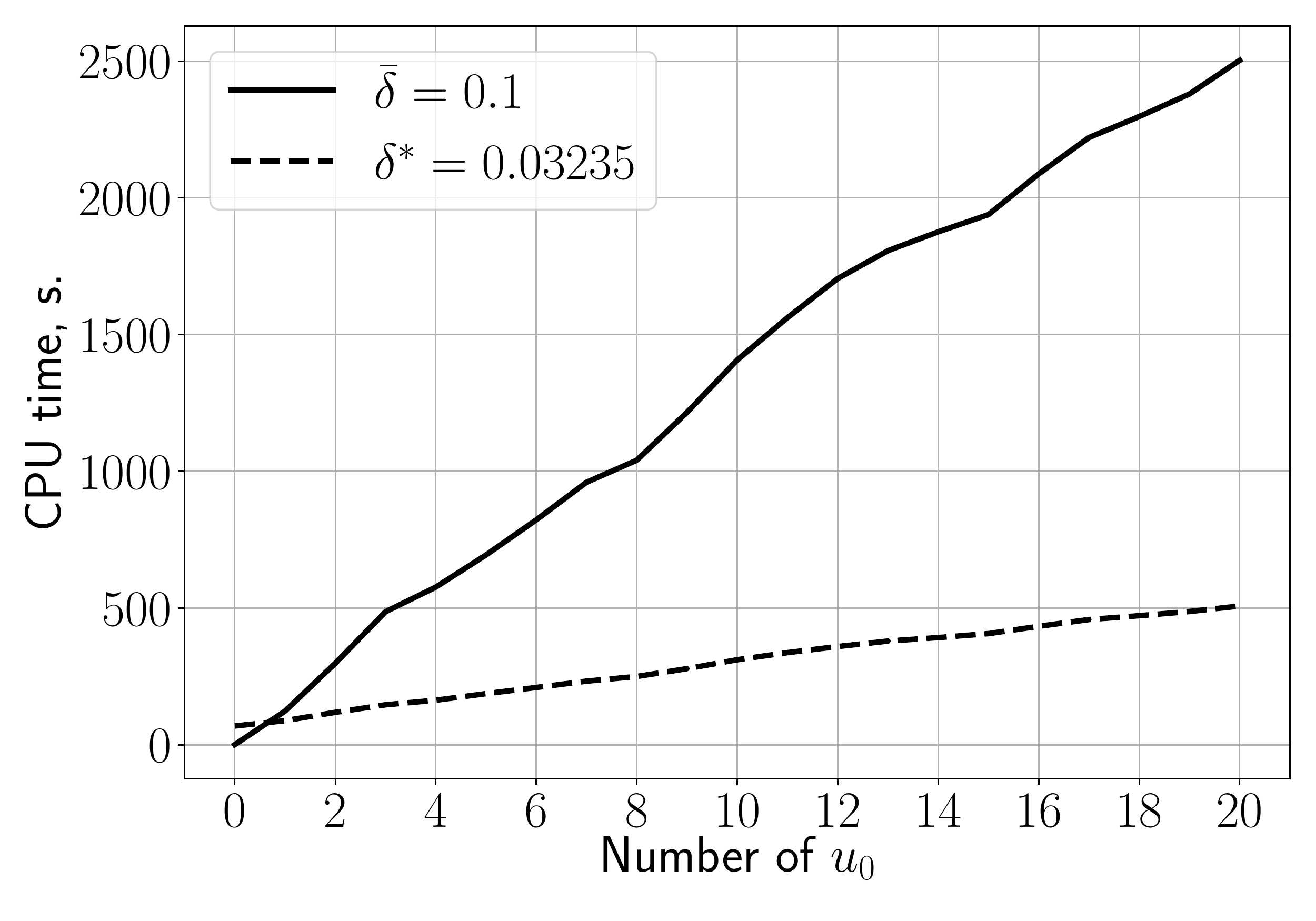}
    \caption{$n = 300, t = 2 \cdot 10^{-4}$}
    \label{fig::opt-run_n300_pe1000_tau2e-4}
    \end{subfigure}
    \caption{Comparison of the total CPU time of the SAI Krylov method with an optimal shift $\gamma^* = \delta^* t$, where $\delta^*$ is determined by the ``optimize-and-run'' method (dashed line) and the SAI Krylov method with a reasonable non-optimized shift $\bar{\gamma} = 0.1t$ (solid line), problem~\eqref{eq::piecewise_diff}. 
    Dashed line starts from the time required to solve~\cref{eq::opt_problem} and find~$\delta^*$.}
    \label{fig::opt-run}
\end{figure}

To demonstrate that the ``optimize-and-run'' method is indeed faster, we provide plots with the total CPU time comparison.
\cref{fig::opt-run} displays how the cumulative CPU time of the SAI Krylov method depends on the number of initial vectors.
The solid lines correspond to the SAI Krylov method with a priori fixed shift $\bar{\gamma}=0.1t$.
The dashed lines correspond to the SAI Krylov method with the optimum shift $\gamma^*$ provided by the ``optimize-and-run'' method.
Note that the starting point of every dashed line is the time required to find $\gamma^*$ for considered~$n$ and~$t$, whereas the solid lines start at the origin.
These plots demonstrate that the optimum shift $\gamma^* = \delta^* t$ provided by the ``optimize-and-run'' method speeds up the convergence of the SAI Krylov method so that the additional costs required to find $\delta^*$ are paid off.
Moreover, we get a gain in the total CPU time already for 4 initial vectors, i.e., $M_{\min} \approx 4$ for this particular problem.
\cref{fig::opt-run} proves our assumption that the costs to find optimum shift $\gamma^* = \delta^*t$ are paid off by the convergence speed up of the SAI Krylov method. 
Also, note that the larger $n$ and $t$, the smaller number of initial vectors are necessary to get a gain in the total CPU time. 
Hence, the ``optimize-and-run'' method is particularly efficient for the fine meshes. 

\paragraph{Effect of $K$ on the total CPU time gain}

\begin{figure}[!ht]
    \centering
    \includegraphics[scale=0.3]{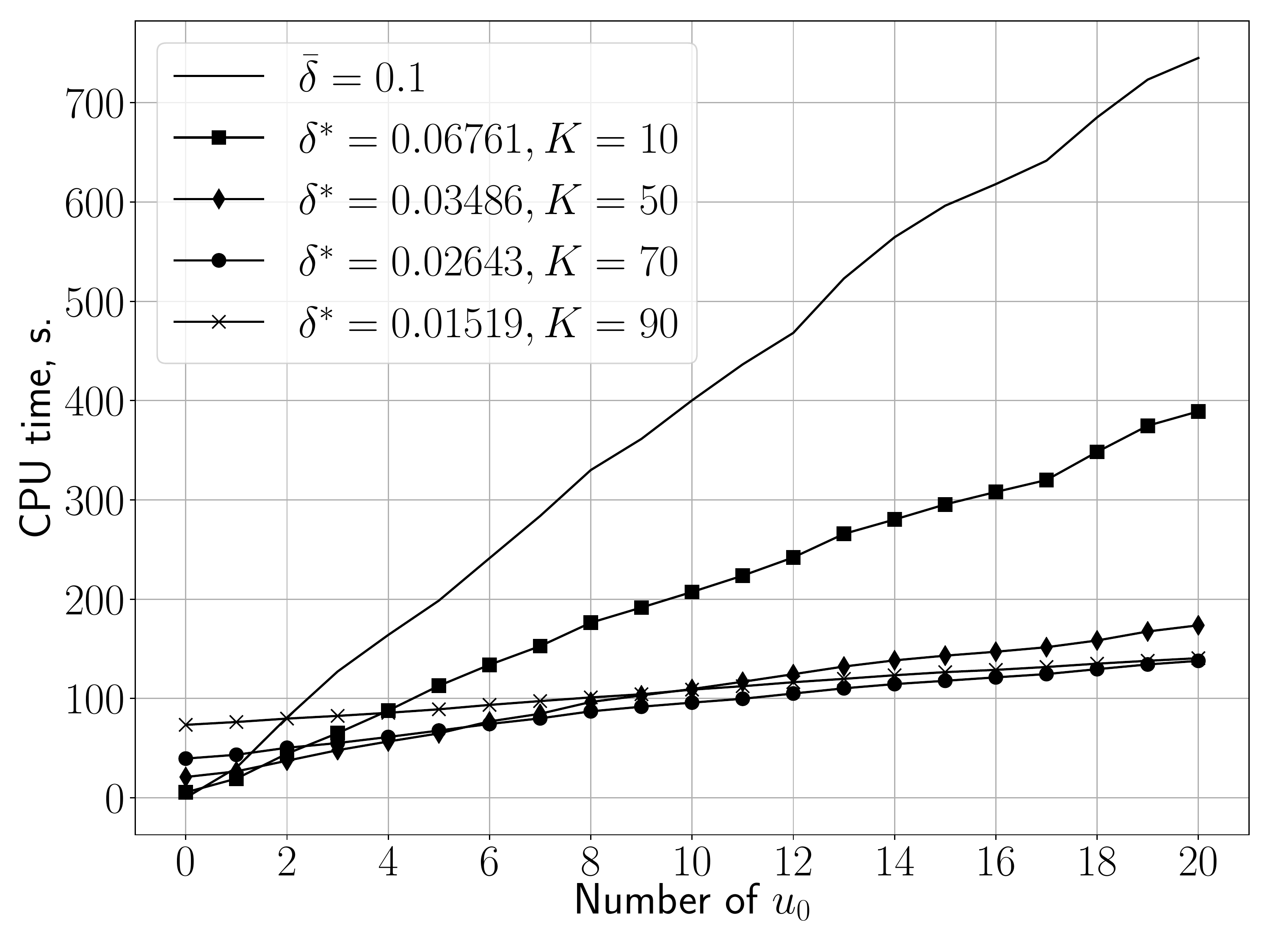}
    \caption{Comparison of the total CPU time of the SAI Krylov method with an optimal shift $\gamma^* = \delta^* t$, where $\delta^*$ determined by the ``optimize-and-run'' method for different $K$ (different marker lines) and with the non-optimized shift $\bar{\gamma} = 0.1 t$ (solid line), problem~\eqref{eq::piecewise_diff}. 
    Matrix $A \in \mathbb{R}^{n^2 \times n^2}$, where $n=200$, and time $t = 4 \cdot 10^{-4}$.}
    \label{fig::K_study}
\end{figure}

As we mentioned above, one of the parameters of the ``optimize-and-run'' method is the number of Arnoldi iterations $K$ used to compute the residual norm for every trial vector~\eqref{eq::opt_problem}.
In~\cref{sec::opt_run} we discuss how $K$ affects the optimization costs, and here we illustrate this numerically.
We consider the case $n=200, t=4\cdot 10^{-4}$ and show the results for four values of $K$: $10, 50, 70, 90$.
For every $K$ we solve  problem~\eqref{eq::opt_problem} and get a  corresponding $\delta^*$.
Then, we measure the total CPU time of processing the initial vectors for every $\delta^*$.
\cref{fig::K_study} presents the comparison of the total CPU time for different  $K$.
This comparison indicates that for larger $K$ we get higher pre-processing costs, but asymptotically faster convergence for a large number of initial vectors.
For example, if $K=10$, then the optimization costs are negligibly small, but the convergence of the SAI Krylov method is not much faster than for a non-optimized shift $\bar{\gamma} = 0.1t$.
Such behavior indicates that the objective function in~\eqref{eq::opt_problem} for small $K=10$ does not properly reflect the convergence of the SAI Krylov method yet.
Hence, the corresponding $\delta^*$ does not guarantee a fast convergence.
At the same time, $K=90$ requires the highest optimization costs, but the gain is almost the same as for $K = 70$.
This example illustrates the importance of a proper choice of~$K$ and a trade-off between optimization costs, obtained convergence speed up and the number of initial vectors to be processed. 

\subsubsection{Incremental method}
\label{sec::increm_piecewise}
\begin{figure}[!htb]
    \centering
    \begin{subfigure}[t]{0.48\textwidth}
    \centering
    \includegraphics[scale=0.27]{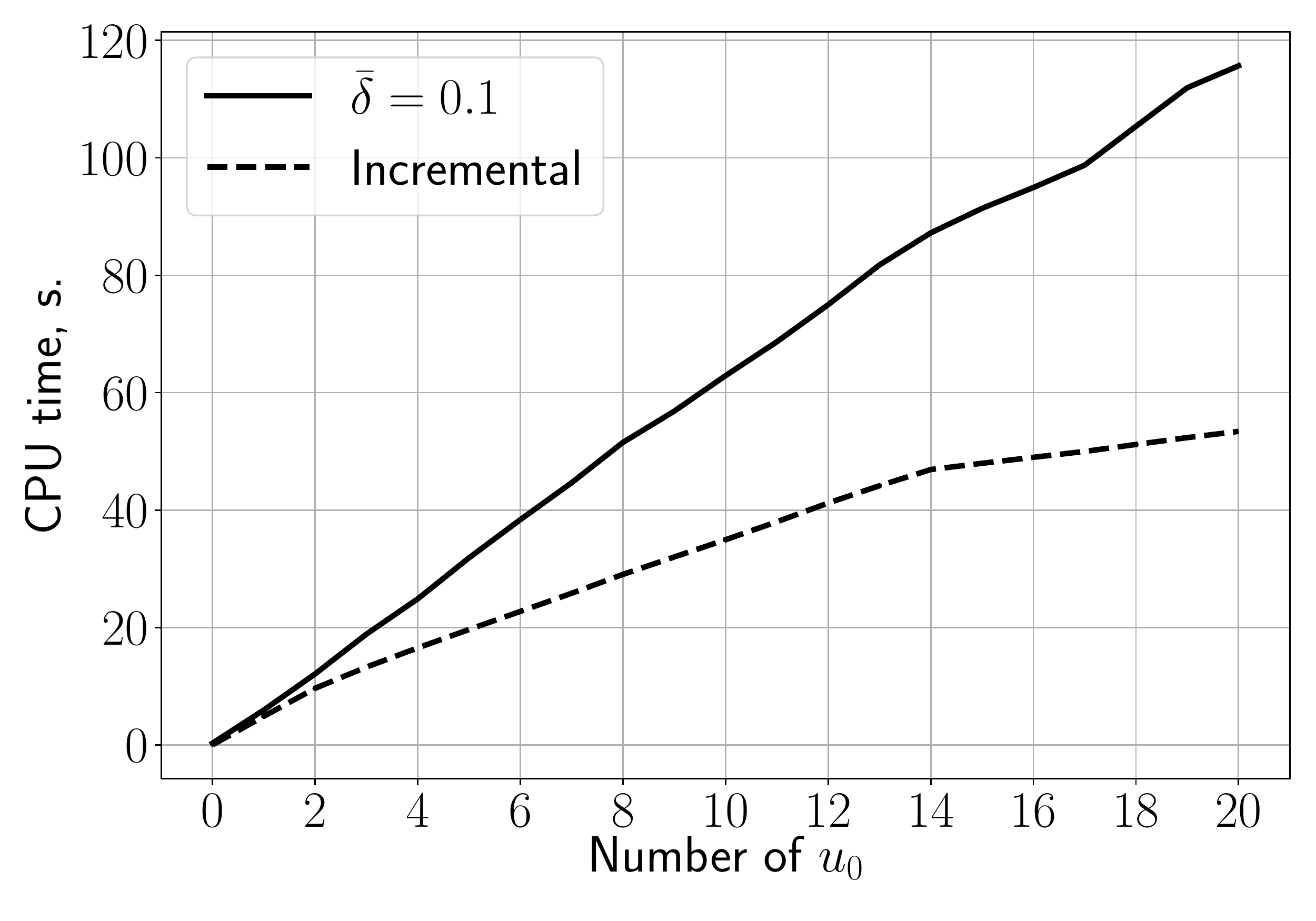}
    \caption{$n = 200, t = 10^{-4}$}
    \label{fig::incremental_n200_pe500_tau0.0001}
    \end{subfigure}
    ~
    \begin{subfigure}[t]{0.48\textwidth}
    \centering
    \includegraphics[scale=0.27]{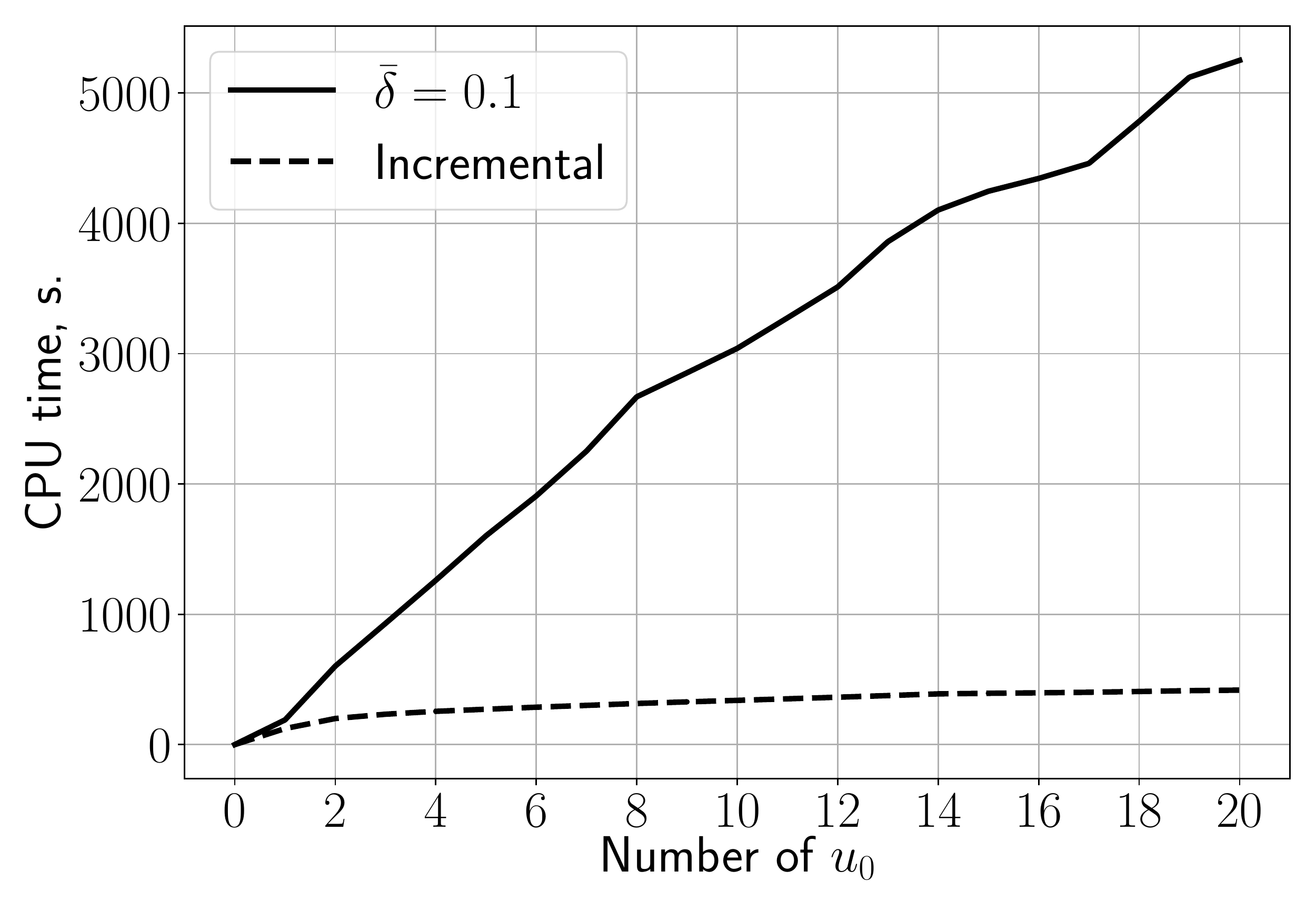}
    \caption{$n=200, t = 4 \cdot 10^{-4}$}
    \end{subfigure}
    \\
    \begin{subfigure}[t]{0.48\textwidth}
    \centering
    \includegraphics[scale=0.27]{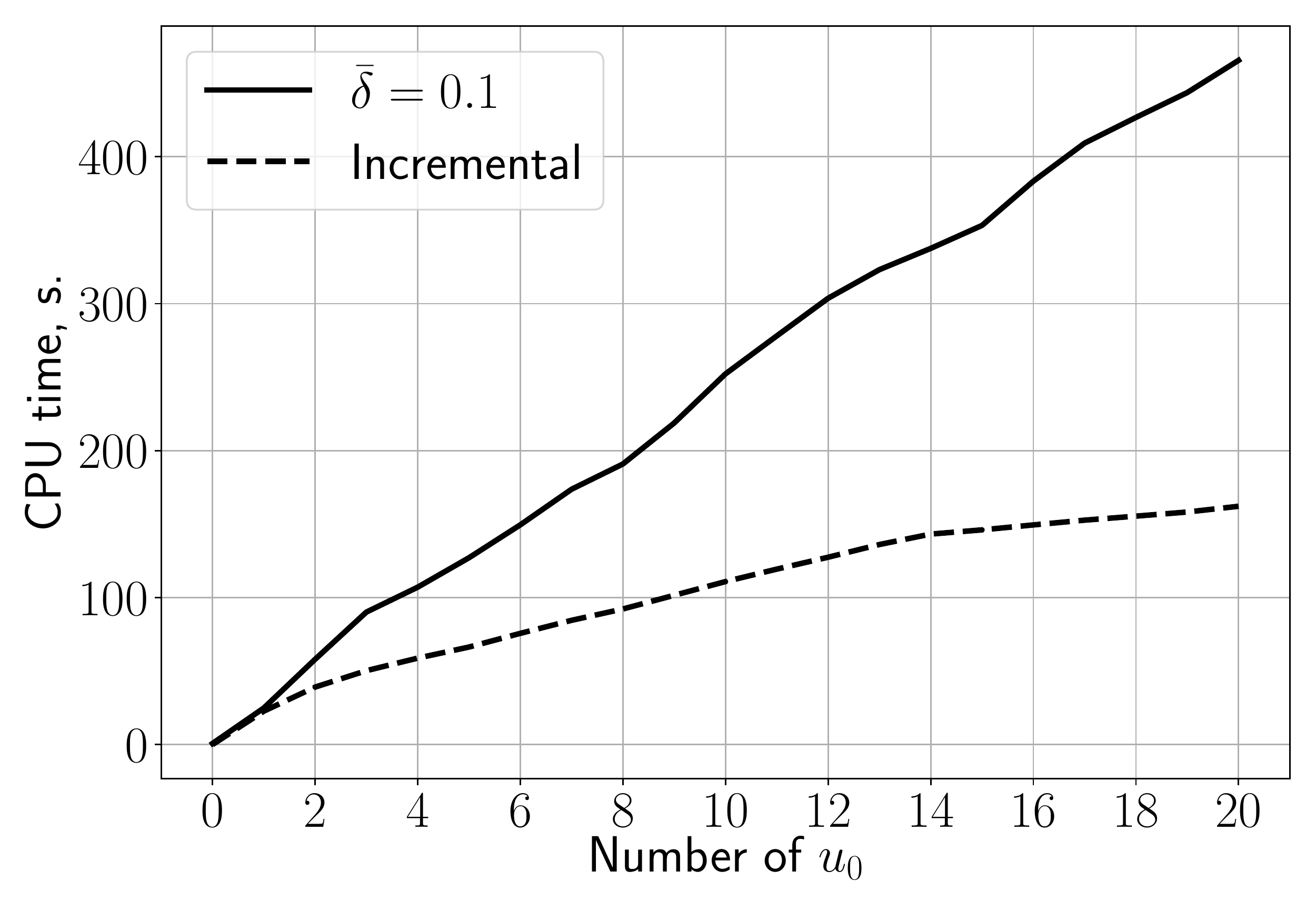}
    \caption{$n=300, t = 10^{-4}$}
    \end{subfigure}
    ~
    \begin{subfigure}[t]{0.48\textwidth}
    \centering
    \includegraphics[scale=0.27]{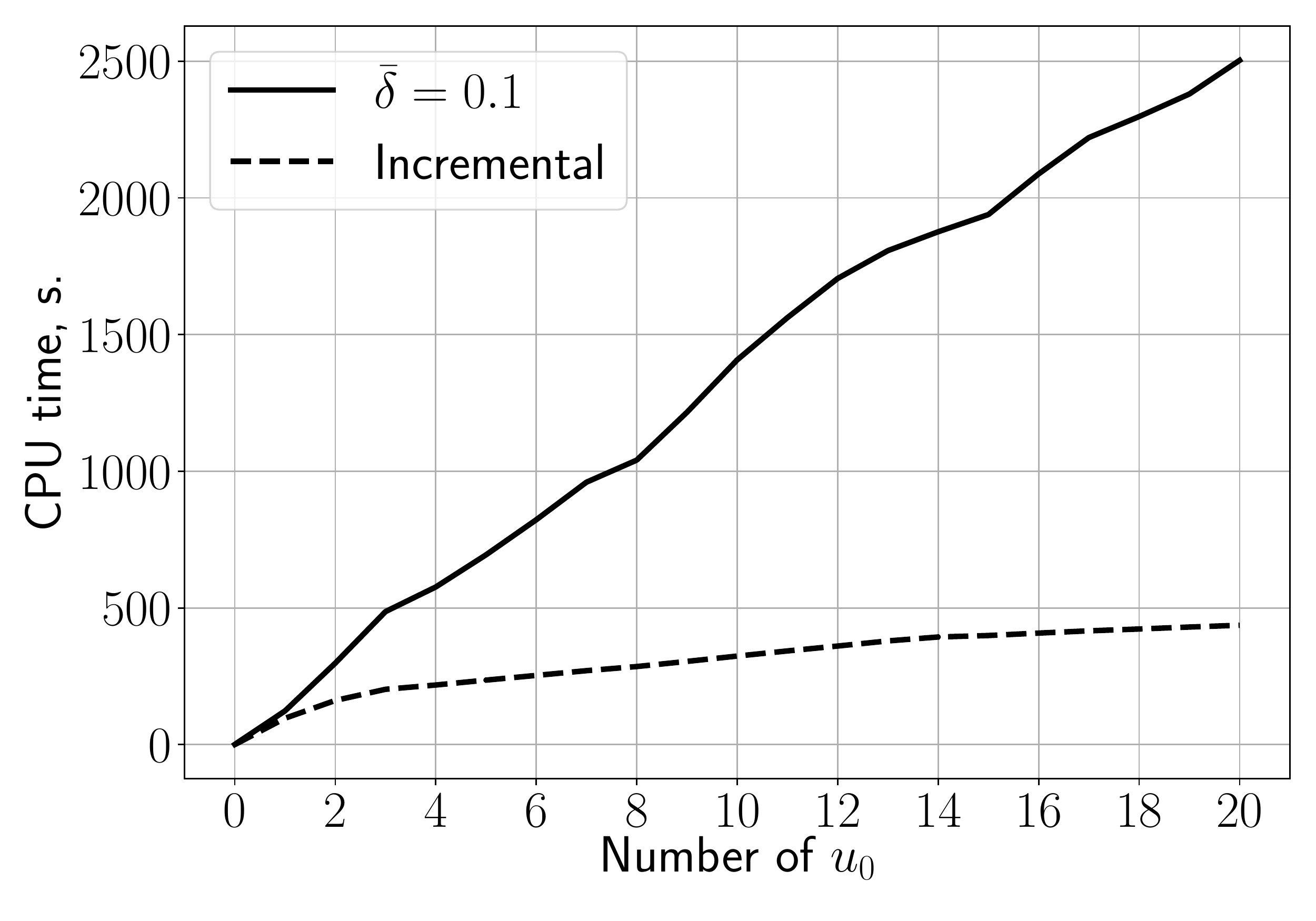}
    \caption{$n=300, t = 2 \cdot 10^{-4}$}
    \end{subfigure}
    \caption{Total CPU time comparison of the SAI Krylov method used incremental tuning of the shift (dashed line) and the SAI Krylov method used constant shift $\bar{\gamma} = 0.1t$ (solid line), problem~\eqref{eq::piecewise_diff}}
    \label{fig::incremental}
\end{figure}

To test the incremental method, we use the same experimental setting as for the ``optimize-and-run'' approach.
\cref{fig::incremental} compares the CPU times of the SAI Krylov method with the constant shift $\bar{\gamma} = 0.1t$ and with the shift values produced by the incremental method.
The same 
% number $N=20$ 
initial vectors are used as in the experiments with the ``optimize-and-run'' method.
Since the incremental method continuously updates the shift value, we do not provide the shift value in the plot.
% \cref{fig::incremental_n200_pe500_tau0.0001} shows that for a small problem size and small $t$ we do not get any significant gain.
% It means that the gain that we get with updated shift values is almost equal to the additional costs to estimate the derivative of the residual norm.
\cref{fig::incremental} demonstrates that we get a significant gain in the CPU time.
Also, the larger $n$ and $t$, the more significant gain we get.
Thus, the costs for shift optimization are paid off and optimized shift values provide sufficient reduction of Arnoldi iteration numbers.     

\subsection{Non-stationary anisotropic diffusion equation}
The second problem we consider is a  non-stationary anisotropic diffusion equation.
This equation with homogeneous Dirichlet boundary conditions and initial conditions is written as
\begin{equation}
    \begin{aligned}
        & \frac{\partial u}{\partial t} = -\mathrm{div}(Q^{\top}\Lambda Q \mathrm{grad} \, u) \quad (x,y)\in\Omega=[-1,1]\times[-1,1] \\
        & u(x, y, 0) = u_0, \quad  u|_{\partial \Omega} = 0,
    \end{aligned}
    \label{eq::aniso}
\end{equation}
where $Q = \begin{bmatrix}
\cos \theta & -\sin \theta \\ \sin \theta & \cos \theta
\end{bmatrix}$ is a rotation matrix and $\Lambda = \begin{bmatrix}
1 & 0 \\ 0 & \lambda 
\end{bmatrix}$ is a matrix representing anisotropy.
To get the matrix $A$ corresponding to the discretization of the right-hand side, we use the second-order finite difference method implemented in the PyAMG library~\cite{OlSc2018}. 
We consider a highly anisotropic case $\lambda = 5000$ and $\theta = \pi/4$.
Note that the matrix $A$ is now symmetric, in contrast to the problem discussed in~\cref{sec::diff_piecewise}. 
Thus, we demonstrate the performance of the proposed methods in another setting.
Since the proposed methods are of interest in the case of multiple initial vectors, we generate 300 initial vectors $v^{(i)}, i=1,\ldots,300$ from the multivariate Gaussian distribution with the same parameters as in~\cref{sec::diff_piecewise}.
The tolerance of the SAI Krylov method is $10^{-8}$.
According to~\cite{van2006preconditioning}, a reasonable shift value for this tolerance is $\bar{\gamma} = \bar{\delta}t = 0.07t$. 
Therefore, $\bar{\delta} = 0.07$, and the search interval in both proposed optimization methods is set to $[0.01, 0.07]$.

\subsubsection{``Optimize-and-run'' method}

In this section performance of the  ``optimize-and-run'' method is compared for problem~\eqref{eq::aniso} with the SAI Krylov method run with reasonable shift $\bar{\gamma}$, which equals $0.07t$ in this case.
Similarly to~\cref{sec::opt_run_discont_diff} we use a single trial vector ($N=1$) to compute the optimum shift $\gamma^* = \delta^*t$ and tolerance in Brent method is set to~$10^{-5}$.
Also, we first discuss the optimization costs and then describe the performance of the ``optimize-and-run'' method.

\paragraph{Optimization cost analysis}
In this section we provide the same optimization cost analysis for the second test problem.
\cref{tab::opt_run_par_n128_aniso,tab::opt_run_par_n256_aniso} provide the total number of sparse $LU$ factorizations, number of Arnoldi iterations $K$ to evaluate the objective function in~\eqref{eq::opt_problem} and the total number of Arnoldi iterations required to find an optimal value $\delta^*$ for $n=128$ and $n=256$, respectively.
These tables show that the number of sparse $LU$ factorizations and Arnoldi iterations are moderate.
Therefore, the optimization costs can be paied off by the obtained convergence speed up.
Also, we note that the larger $n$ does not require the larger number of sparse $LU$ factorizations $s$ or number of Arnoldi iterations $K$ to evaluate the objective function in~\eqref{eq::opt_problem}.
This is expected because the SAI Krylov method should, in principle, exhibit a mesh independent convergence.
In experiments we use $K=20$ and $s\approx 16$ for $n=128$ and $n=256$.

\begin{table}[!htb]
    \centering
    \caption{Costs of the ``optimize-and-run'' method to find a proper shift value $\gamma^* = \delta^*t$ for $n = 128$, $A \in \mathbb{R}^{n^2 \times n^2}$, problem~\eqref{eq::aniso}}
    \begin{tabular}{ccc}
    \toprule
        & $t = 0.1$ & $t = 0.5$ \\
        \midrule
    Total number of sparse $LU$ factorizations, $s$ &  17  &  18 \\
    Number of Arnoldi iterations for every initial vector, $K$ & 20 & 20 \\
    Total number of Arnoldi iterations  & $ 17\cdot 20$ & $ 18\cdot 20$ \\
    \bottomrule
    \end{tabular}
    \label{tab::opt_run_par_n128_aniso}
\end{table}

\begin{table}[!htb]
    \centering
    \caption{Costs of the ``optimize-and-run'' method to find a proper shift value $\gamma^* = \delta^*t$ for $n = 256$, $A \in \mathbb{R}^{n^2 \times n^2}$, problem~\eqref{eq::aniso}}
    \begin{tabular}{ccc}
    \toprule
        & $t = 0.1$ & $t = 0.5$  \\
        \midrule
    Total number of sparse $LU$ factorizations, $s$ &  $15$ &  $17$ \\
    Number of Arnoldi iterations for every initial vector, $K$ & $20$ & $20$  \\
    Total number of Arnoldi iterations & $15\cdot20$ & $17\cdot20$ \\
    \bottomrule
    \end{tabular}
    \label{tab::opt_run_par_n256_aniso}
\end{table}

\paragraph{Performance comparison}
Now we show that the ``optimize-and-run'' method provides a better shift value for the SAI Krylov method than an a priori chosen reasonable shift value~$\bar{\gamma}$.
To do so, we measure the average number of Arnoldi iterations per initial vector and total CPU time to process the all initial vectors. 
\cref{tab::opt_run_arnoldi_test_aniso} presents the average number of Arnoldi iterations in the SAI Krylov method with these shift values.
As we see, the optimum shift $\gamma^* = \delta^*t$ indeed leads to a reduction in the average number of Arnoldi iterations.
This indicates that asymptotically, for a growing number of initial vectors, we get a gain in the total CPU time, too.
The minimum number of initial vectors $M_{\min}$ to observe this gain is provided in the next paragraph.     
However, the observed gain for this problem is much smaller than we get for the problem~\cref{eq::piecewise_diff}, cf.~\cref{tab::opt_run_arnoldi_test}.
It means that the dependence of the number of Arnoldi iterations on the value shift is not so strong as it is in the problem~\eqref{eq::piecewise_diff}.
A possible explanation of this weak dependence is that the problem is symmetric in contrast to the problem~\eqref{eq::piecewise_diff}.

\begin{table}[!ht]
    \centering
    \caption{Average number of Arnoldi iterations over 300 initial vectors to achieve tolerance $10^{-8}$, matrix $A \in \mathbb{R}^{n^2 \times n^2}$, problem~\eqref{eq::aniso}}
    \begin{tabular}{ccc}
    \toprule
       & $\bar{\delta} = 0.07$ & $\delta^*$  \\
       \midrule
     $n=128, t=0.1$  & $30.72$ & $28.15$ \\
     $n=128, t=0.5$ & $30.96$  & $28.18$ \\
     $n=256, t=0.1$ & $30.44$  & $26.99$  \\
     $n=256, t=0.5$ & $30.63$ & $27.32$  \\
     \bottomrule
    \end{tabular}
    \label{tab::opt_run_arnoldi_test_aniso}
\end{table}

\begin{figure}[!htb]
    \centering
    \begin{subfigure}[t]{0.48\textwidth}
    \centering
    \includegraphics[scale=0.27]{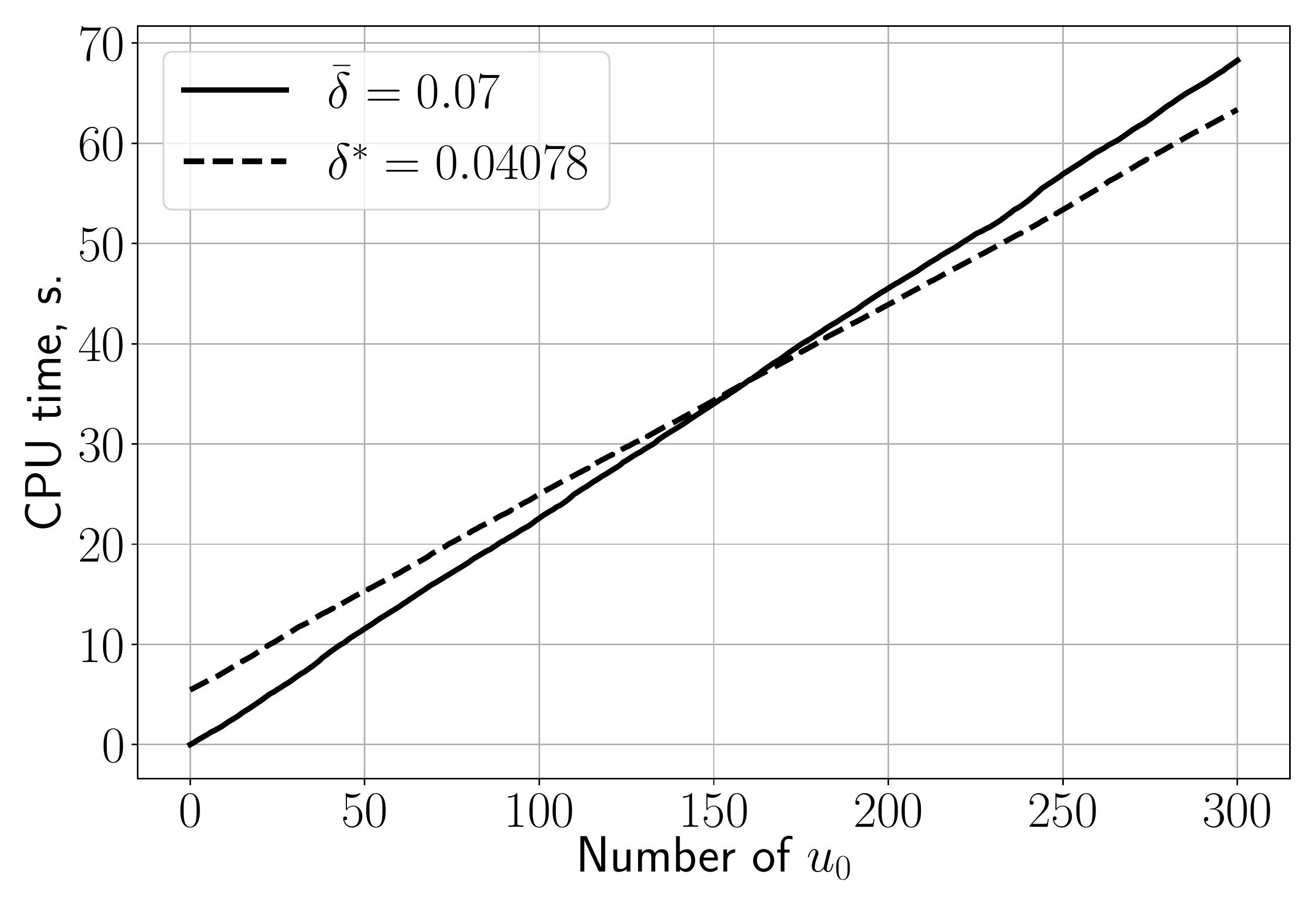}
    \caption{$n = 128, t = 0.1$}
    \label{fig::opt-run_n128_tau0.1_aniso}
    \end{subfigure}
    ~
    \begin{subfigure}[t]{0.48 \textwidth}
    \centering
    \includegraphics[scale=0.27]{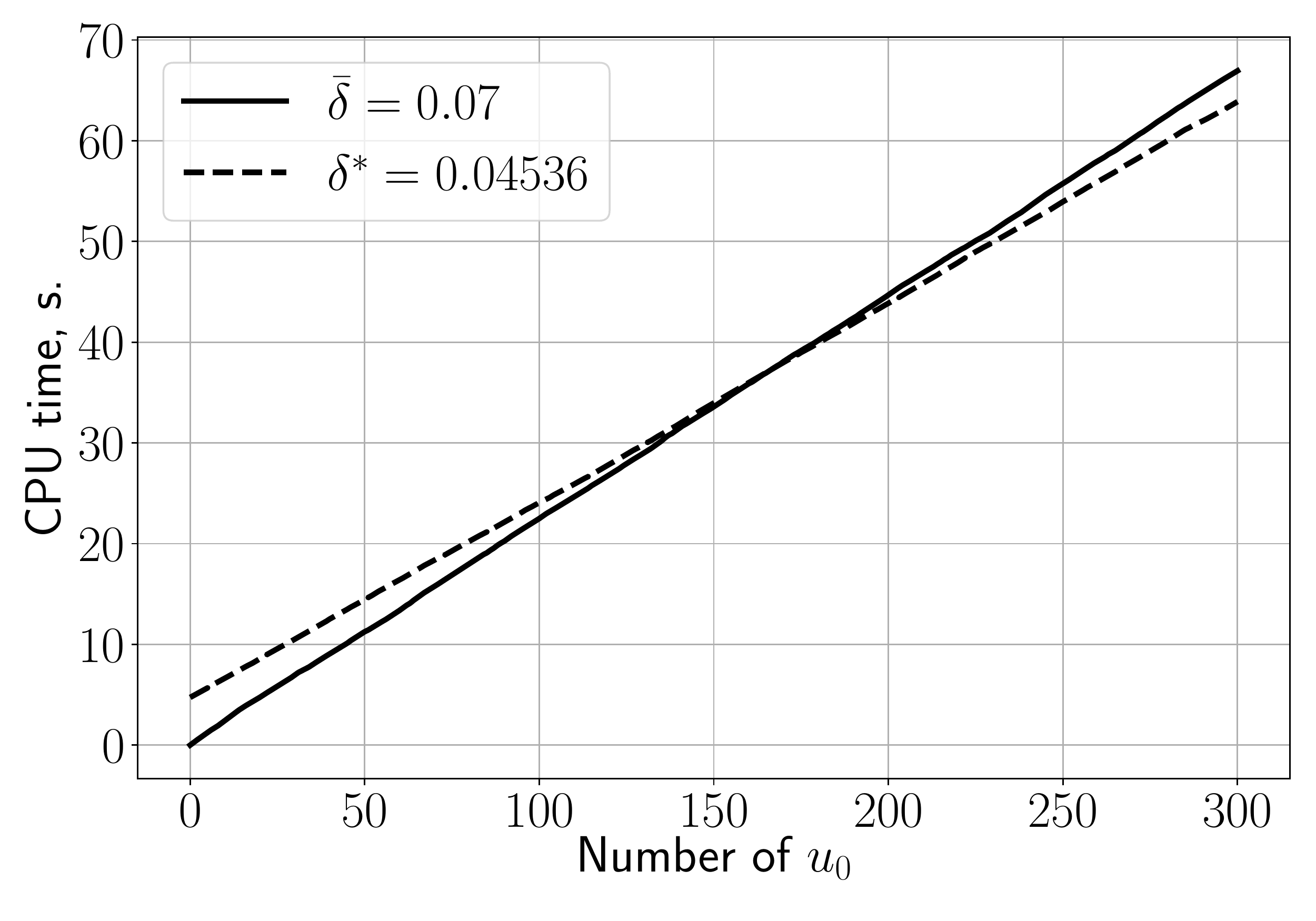}
    \caption{$n = 128, t = 0.5$}
    \end{subfigure}
    \\
    \begin{subfigure}[t]{0.48\textwidth}
    \centering
    \includegraphics[scale=0.27]{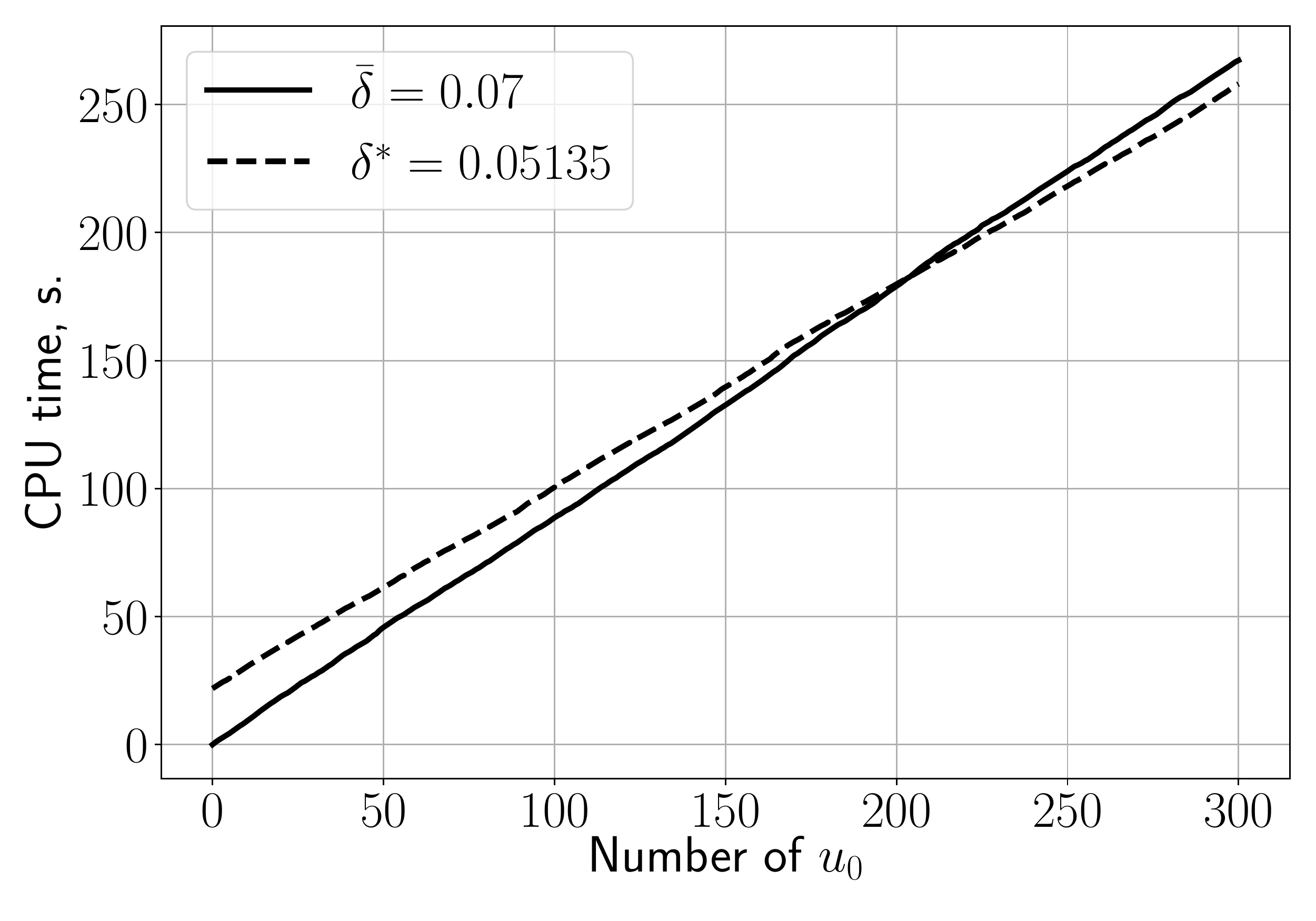}
    \caption{$n = 256, t = 0.1$}
    \label{fig::opt-run_n256_tau0.1_aniso}
    \end{subfigure}
    ~
    \begin{subfigure}[t]{0.48\textwidth}
    \centering
    \includegraphics[scale=0.27]{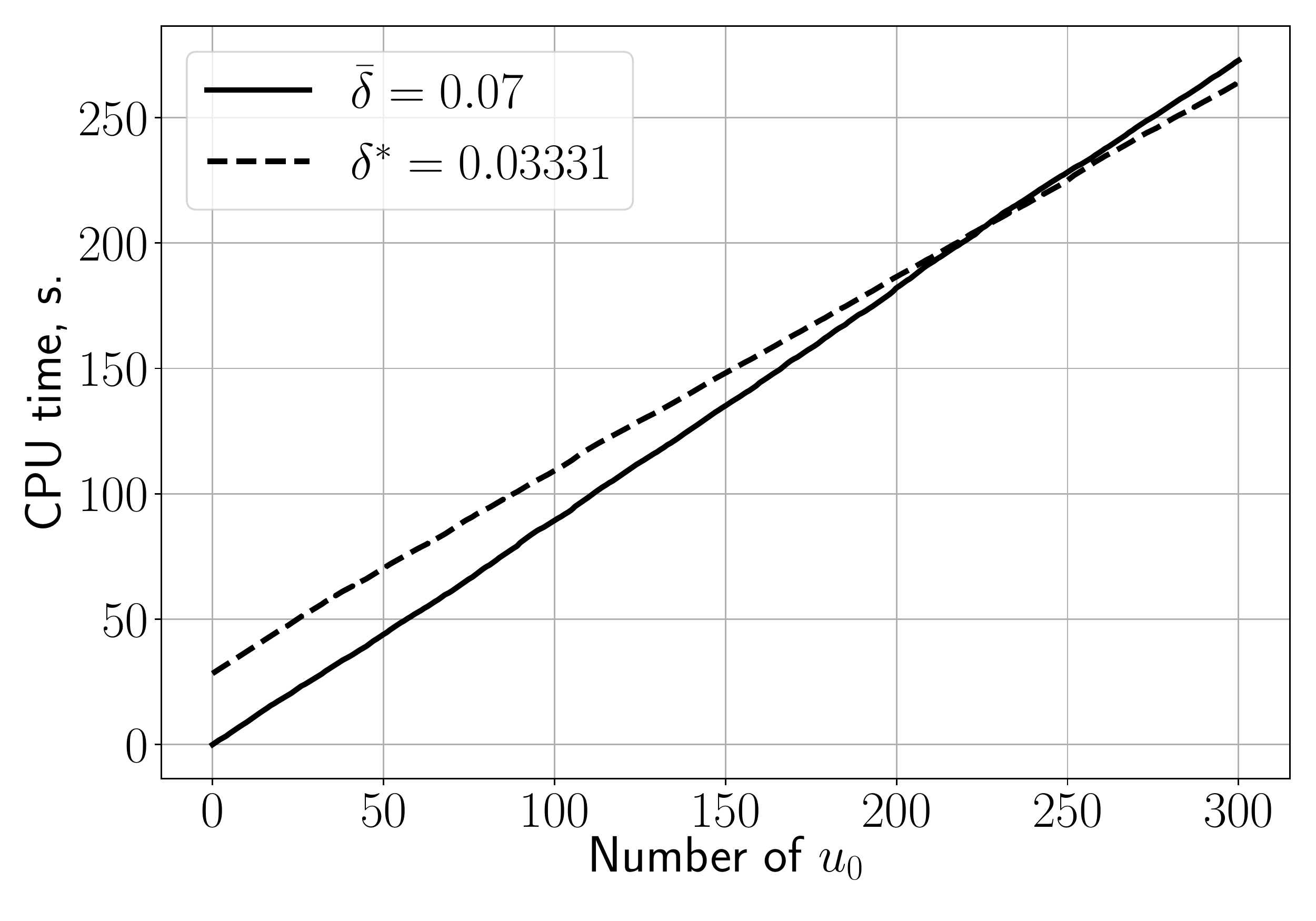}
    \caption{$n = 256, t = 0.5$}
    \label{fig::opt-run_n256_tau0.5_aniso}
    \end{subfigure}
    % \\
    % \begin{subfigure}[t]{0.45\textwidth}
    % \centering
    % \includegraphics[scale=0.23]{{tau0.1_batch1_iter20_aniso5000_theta45_n512_tol1e-8}.pdf}
    % \caption{$n = 512, t = 0.1$}
    % \end{subfigure}
    % ~
    % \begin{subfigure}[t]{0.45\textwidth}
    % \centering
    % \includegraphics[scale=0.23]{{tau0.5_batch1_iter20_aniso5000_theta45_n512_tol1e-8}.pdf}
    % \caption{$n = 512, t = 0.5$}
    % \end{subfigure}
    \caption{Comparison of the total CPU time of the SAI Krylov method with an optimal shift $\gamma^* = \delta^* t$, where $\delta^*$ is determined by the ``optimize-and-run'' method (dashed line) and the SAI Krylov method with a reasonable non-optimized shift $\bar{\gamma} = 0.07t$ (solid line), problem~\eqref{eq::aniso}. 
    Dashed line starts from the time required to solve~\cref{eq::opt_problem} and find~$\delta^*$.}
    \label{fig::opt-run_aniso}
\end{figure}

\cref{fig::opt-run_aniso} shows the total CPU time of the SAI Krylov method with the considered shifts.
The optimum shift $\gamma^*$ given by the ``optimize-and-run'' method provides a faster processing of the initial vectors in all considered experimental settings.
Thus, the optimization costs are paid off. 
Also, from this plot, we see that the approximate minimum number of initial vectors $M_{\min}$ to get a gain in the total CPU time varies from 170 to 210.
This value of $M_{\min}$ is much bigger than for the problem~\eqref{eq::piecewise_diff}.
This agrees well with the observation on the average number of Arnoldi iterations reduction, see~\cref{tab::opt_run_arnoldi_test_aniso}.

\subsubsection{Incremental method}
\begin{figure}[!htb]
    \centering
    \begin{subfigure}[t]{0.48\textwidth}
    \centering
    \includegraphics[scale=0.27]{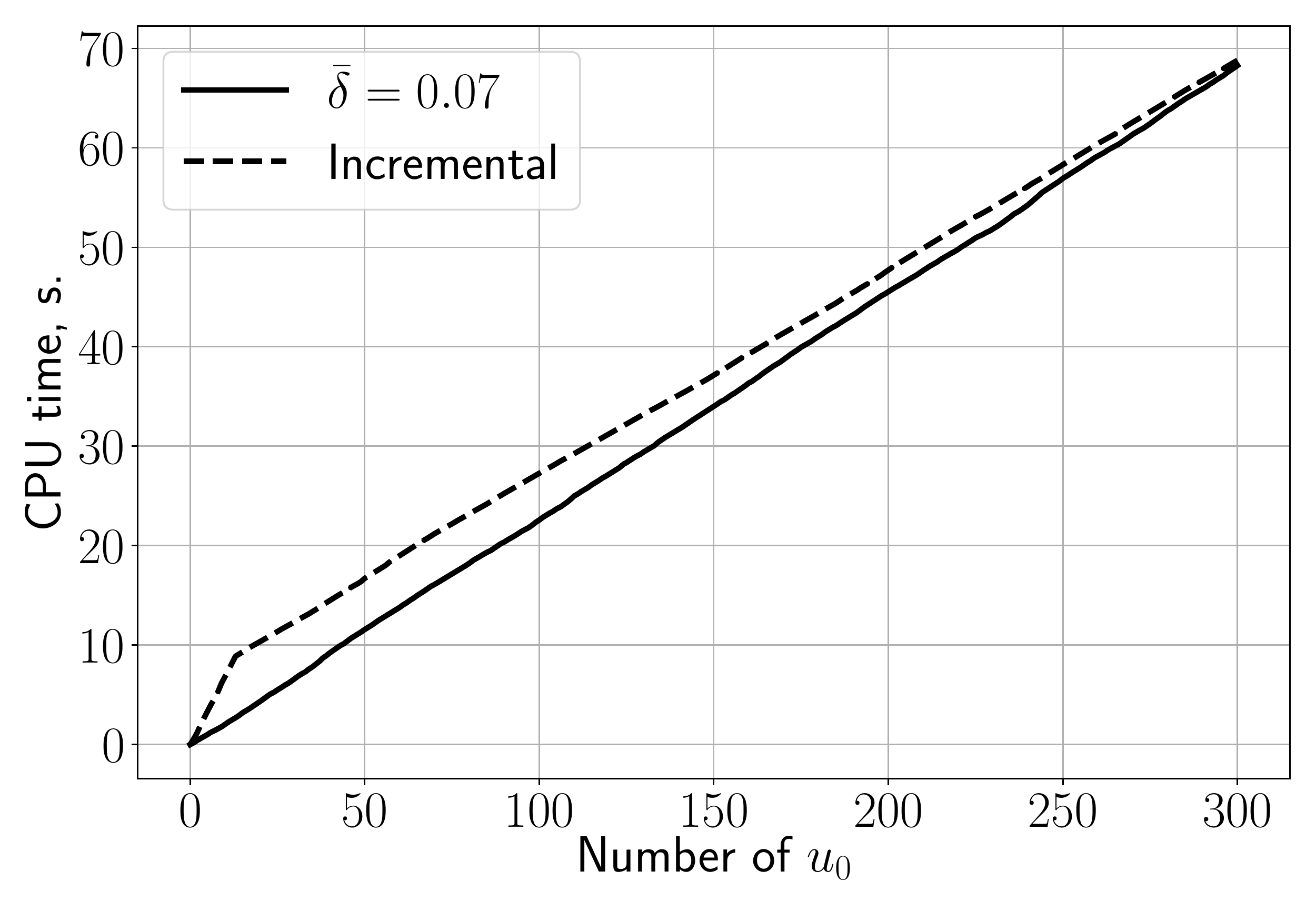}
    \caption{$n = 128, t = 0.1$}
    \label{fig::incremental_n128_tau0.1}
    \end{subfigure}
    ~
    \begin{subfigure}[t]{0.48\textwidth}
    \centering
    \includegraphics[scale=0.27]{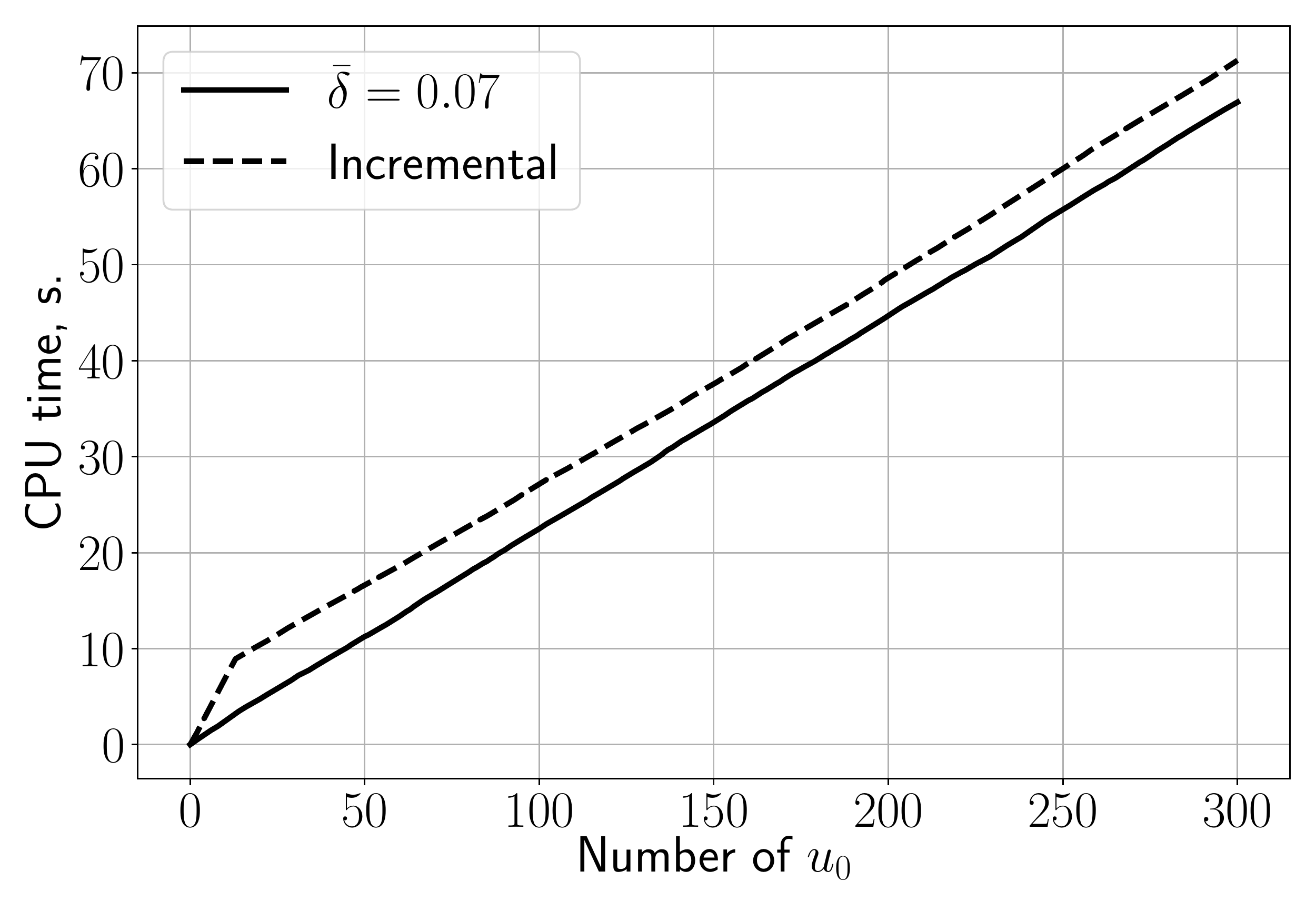}
    \caption{$n=128, t = 0.5$}
    \label{fig::incremental_n128_tau0.5}
    \end{subfigure}
    \\
    \begin{subfigure}[t]{0.48\textwidth}
    \centering
\includegraphics[scale=0.27]{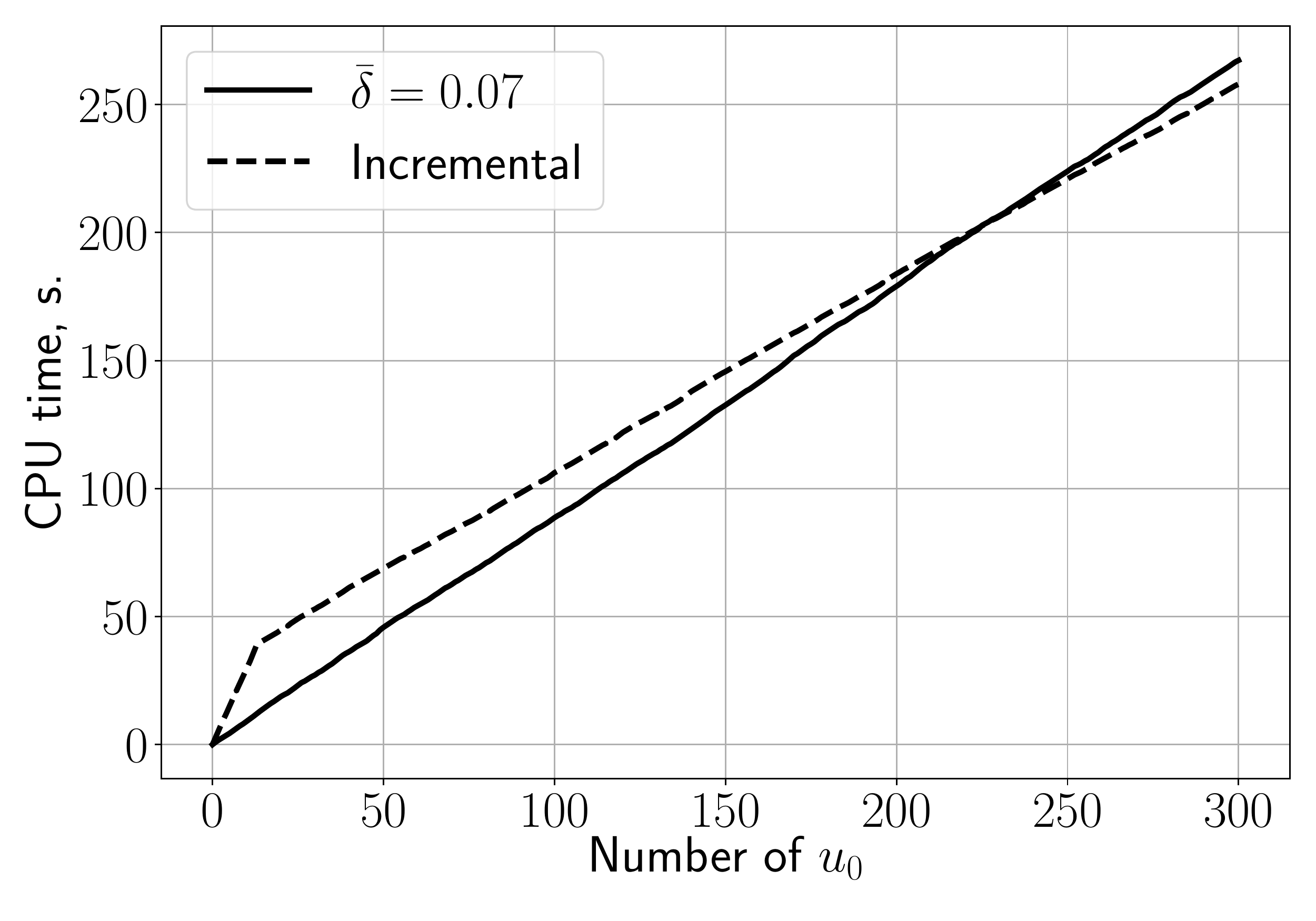}
    \caption{$n=256, t = 0.1$}
    \label{fig::incremental_n256_tau0.1}
    \end{subfigure}
    ~
    \begin{subfigure}[t]{0.48\textwidth}
    \centering
    \includegraphics[scale=0.27]{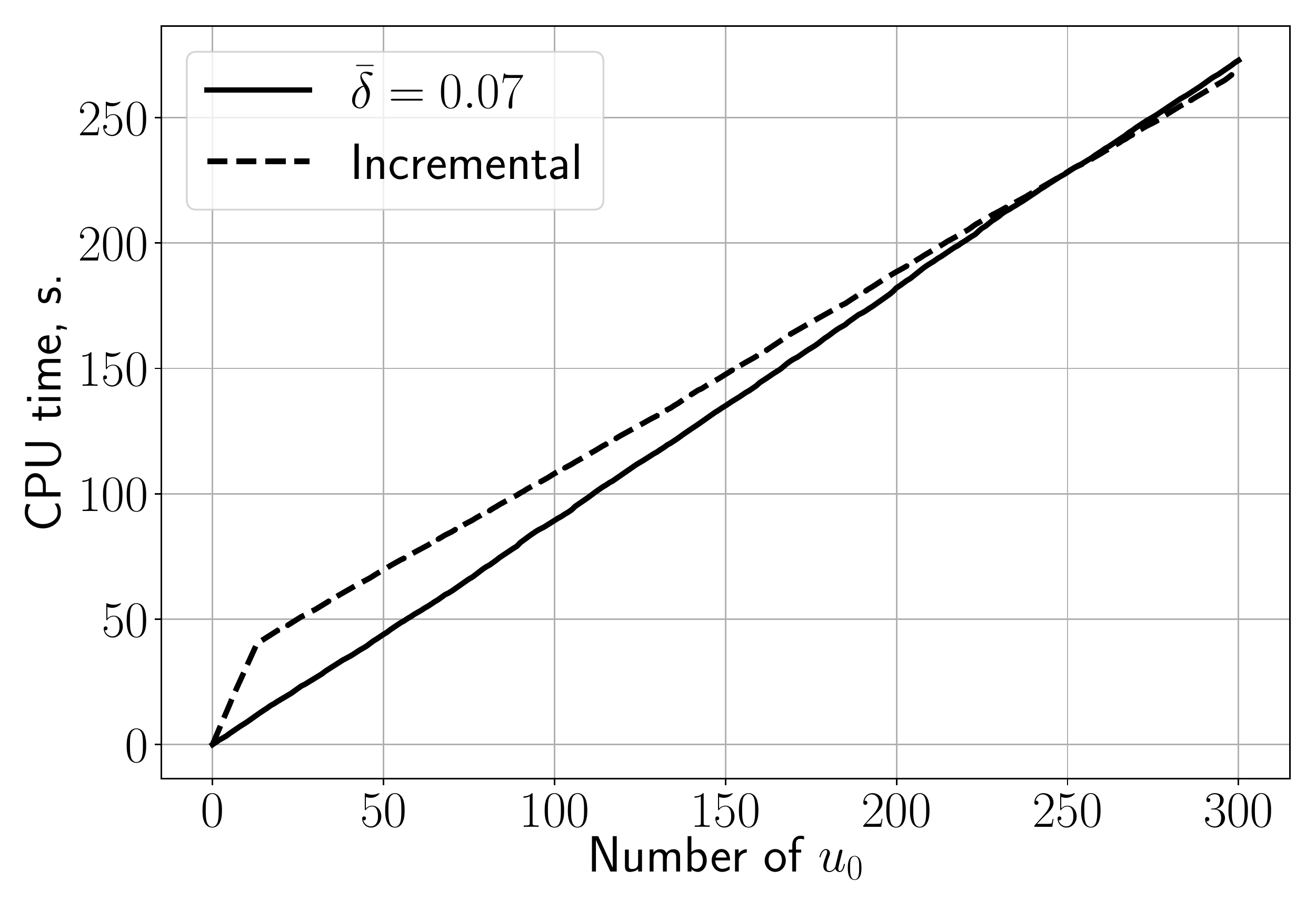}
    \caption{$n=256, t = 0.5$}
    \label{fig::incremental_n256_tau0.5}
    \end{subfigure}
    \caption{Total CPU time comparison of the SAI Krylov method used incremental tuning of the shift (dashed line) and the SAI Krylov method used constant shift $\bar{\gamma} = 0.07t$ (solid line), problem~\eqref{eq::aniso}}
    \label{fig::incremental_aniso}
\end{figure}

In this section, we test the incremental method for the considered anisotropic diffusion equation~\eqref{eq::aniso}.
% To reduce the computational costs, we slightly modify the incremental method presented in~\cref{alg::update_delta}.
% The modified incremental method stops updating  $\delta$ if the difference between two sequential $\delta$ is less in modulus than $10^{-5}$.
% This modification reduces the number of sparse $LU$ factorizations and, consequently, the total running time.
% So, we have two stages in the modified incremental method. 
% At the first stage the method proceses initial vectors and updates $\delta$ according to~\cref{alg::update_delta} until it has converged to some value which we denote $\tilde{\delta}$. 
% At the second stage the SAI Krylov method runs with the constant shift $\tilde{\gamma} = \tilde{\delta}t$ for the remaining initial vectors without derivative estimation.
\cref{fig::incremental_aniso} presents the performance of the SAI Krylov method with the incremental shift tuning and with a constant shift $\bar{\gamma} = 0.07t$.
This plot shows the dependence of the total CPU time of the SAI Krylov method on the number of initial vectors.
% The solid lines correspond to the constant shift $\bar{\gamma}$ and the dashed lines correspond to the incremental tuning of shift. 
\cref{fig::incremental_n256_tau0.1,fig::incremental_n256_tau0.5} demonstrate that the incremental shift tuning leads to faster processing of 300 or more initial vectors than the constant shift $\bar{\gamma} = 0.07t$.
Similarly to~\cref{sec::increm_piecewise}, the larger the number of initial vectors, the more significant gain we get.
At the same time, \cref{fig::incremental_n128_tau0.1} shows that, for setting $n=128$ and $t=0.1$, the incremental shift tuning leads to a processing speed for more than 400 initial vectors.
% In particular, the slope of the dashed line is smaller than the one for the solid line, but 300 initial states are not enough to get the gain in the total running time.
For the setting $n=128$ and $t=0.5$,~\cref{fig::incremental_n128_tau0.5}, the incremental shift tuning does not give a convergence speed up of the SAI Krylov method. 
The dashed and solid lines are parallel, and therefore we do not get even an asymptotic gain.
A possible explanation of this result is that the derivative estimation for given first initial states is not accurate enough.

% The second order finite difference scheme gives the following stencil
% \[
% \begin{bmatrix}
% 1/2(\lambda - 1) \cos \theta \sin \theta &    -(\cos^2 \theta + \lambda \sin^2\theta) &     -1/2(\lambda - 1) \cos \theta \sin \theta\\ 
% -(\lambda \cos^2 \theta + \sin^2 \theta) &      2 (\lambda + 1) &    -(\lambda \cos^2 \theta + \sin^2 \theta)    \\
% -1/2(\lambda - 1) \cos \theta \sin \theta &   -(\cos^2 \theta + \lambda \sin^2 \theta)  & 1/2(\lambda - 1) \cos \theta \sin \theta 
% \end{bmatrix}
% \]

\section{Conclusions}
\label{sec:conclusion}

In this paper, we consider the problem of a proper choice of the shift in the SAI Krylov method for computing matrix-vector products with the matrix exponential. 
To choose the shift, we propose the ``optimize-and-run'' method and the incremental method. 
These methods prove to be useful if the products with the matrix exponential have to be computed for a number of vectors. 
% MB:
In the experiments, these initial vectors are taken to be random with normally distributed entries.
The methods are complementary to each other in the sense that they handle these vectors differently.
In particular, the ``optimize-and-run'' method is designed for the case where many vectors are available in advance and we can use or generate them to perform the optimization stage.
The optimization stage requires a few trial vectors, but we observe that only one trial vector is enough to get a sufficiently good optimal shift in the considered types of problems.
The optimal shift found at optimization stage is then used in processing the other vectors.
In contrast, the incremental method does not require that the all the initial vectors are available beforehand and can process the vectors one by one.
To demonstrate the performance of the proposed methods we consider two test problems.
The non-symmetric matrix exponential action on some initial vectors is computed in the first test problem and the symmetric matrix exponential action in the second one.
In both test problems, the proposed methods give a gain in the total CPU time compared to the SAI Krylov method run with a reasonable shift value.
It turns out that in some cases the additional costs in both methods are paid off already for a quite moderate number of initial vectors (for instance, for~3 or 4~vectors).
%MB: added 
These conclusions hold provided all the initial vectors belong to a set of random vectors with normally distributed entries.
Thus, both proposed methods for a proper shift choice require moderate additional costs and produce shift values which provide a faster convergence of the SAI Krylov method than other reasonable shift values.
% To test the proposed methods, we consider the application of the matrix exponential to the integration of the partial differential equations (PDEs).
% In tests, we study two examples of practically important PDEs: the diffusion equation with piecewise constant coefficients and the anisotropic diffusion equation.
% The proposed methods demonstrate better performance for the considered problems than the other approach to choose the shift. 
% We discuss the obtained results and the feature of the test problems that leads to gain in the total running time.
% Possible future work is to extend the proposed methods to the SAI Krylov method with multiple shifts.

\bibliographystyle{siamplain}
\bibliography{lib}
\end{document}